%% file: ms.tex
\documentclass[twoside,english]{article}
\usepackage[T1]{fontenc}
\usepackage[latin9]{inputenc}
\usepackage{color}
\usepackage{babel}
\usepackage{verbatim}
\usepackage{amsmath}
\usepackage{amsthm}
\usepackage{amssymb}
\usepackage{graphicx}
\usepackage{epstopdf}
\ifpdf
\DeclareGraphicsExtensions{.eps,.pdf,.png,.jpg}
\else
\DeclareGraphicsExtensions{.eps}
\fi
\usepackage{geometry}
\usepackage{setspace}
\usepackage{wasysym}
\usepackage[pdfusetitle,
 bookmarks=true,bookmarksnumbered=false,bookmarksopen=false,
 breaklinks=false,pdfborder={0 0 1},backref=false,colorlinks=true]
 {hyperref}
\hypersetup{
 citecolor=cyan}

\makeatletter

\let\SF@@footnote\footnote
\def\footnote{\ifx\protect\@typeset@protect
    \expandafter\SF@@footnote
  \else
    \expandafter\SF@gobble@opt
  \fi
}
\expandafter\def\csname SF@gobble@opt \endcsname{\@ifnextchar[%]
  \SF@gobble@twobracket
  \@gobble
}
\edef\SF@gobble@opt{\noexpand\protect
  \expandafter\noexpand\csname SF@gobble@opt \endcsname}
\def\SF@gobble@twobracket[#1]#2{}
\theoremstyle{plain}
\newtheorem{prop}{\protect\propositionname}
\theoremstyle{plain}
\newtheorem{thm}{\protect\theoremname}
\theoremstyle{remark}
\newtheorem*{rem*}{\protect\remarkname}
\theoremstyle{plain}
\newtheorem{cor}{\protect\corollaryname}
\theoremstyle{plain}
\newtheorem{assumption}{\protect\assumptionname}
\theoremstyle{definition}
\newtheorem{defn}{\protect\definitionname}
\theoremstyle{plain}
\newtheorem{lem}{\protect\lemmaname}
\theoremstyle{remark}
\newtheorem{rem}{\protect\remarkname}

\providecommand{\casualprint}{}%
\usepackage{lmodern}
\usepackage[dvipsnames]{xcolor}
\usepackage{tikz-cd}

\newcommand{\sepia}[1]{#1}

\ifdefined\showcaptionsetup
 \PassOptionsToPackage{caption=false}{subfig}
\fi
\usepackage{subfig}
\makeatother

\providecommand{\assumptionname}{Assumption}
\providecommand{\corollaryname}{Corollary}
\providecommand{\definitionname}{Definition}
\providecommand{\lemmaname}{Lemma}
\providecommand{\propositionname}{Proposition}
\providecommand{\remarkname}{Remark}
\providecommand{\theoremname}{Theorem}

\begin{document}
\include{math_shorthand}

\title{A Dilation-based Seamless Multiscale Method For Elliptic Problems}
\author{Ziheng Chen, Bj\"{o}rn Engquist}
\maketitle
\begin{abstract}
Many numerical methods for multiscale differential equations require a scale separation between the larger and the smaller scales to achieve accuracy and computational efficiency. In the area of multiscale dynamical systems, so-called, seamless methods have been introduced to reduce the requirement of scale separation. We will translate these methods to numerical homogenization problems and extend the technique to multiple dimensions. The initial step is to prove that a one-dimensional \sepia{second-order} elliptic operator with oscillatory coefficients can be rewritten as a multiscale dynamical system. Inspired by this, multiscale elliptic operators in higher dimensions are approximated by a novel approach based on local dilation, which provides a middle ground for balancing intractability and accuracy without the need for full resolution. The dilation operator can be further generalized to preserve important structures by properly decomposing the coefficient field. Error estimates are developed and promising numerical results of different examples are included.

\ifdefined\casualprint{}\else{
}\fi
\end{abstract}

\section{Introduction\label{sec:Introduction}}

Multiscale problems have long been challenging in numerical analysis due to the computational cost of resolving the smallest scale over the range of the largest scales. Although the homogenization theory is able to characterize the behavior on the coarse scale, it usually requires solving a family of PDEs which is often computationally intractable. We are interested in the following \sepia{second-order} elliptic problem
\begin{equation}
-\nabla\cdot\left[A^{\epsilon}\left(x\right)\nabla u^{\epsilon}\left(x\right)\right]=f\left(x\right)\ \text{in}\ \Omega\label{eq:intro-elliptic-pde}
\end{equation}
with the Dirichlet boundary condition. To model the multiscale behavior, we assume the diffusion coefficient takes the form $A^{\epsilon}\left(x\right)=A\left(x,\frac{x}{\epsilon}\right)$ where $A$ is a 1-periodic function in the second argument. By \cite{Bensoussan1978,Allaire1992,Grigorios2008,Suslina2013,Zhikov2016}, the sequence of solutions $u^{\epsilon}$ converges, in the weak $H^{1}$ topology as $\epsilon\to0$, to the solution $u$ of the homogenized equation
\begin{equation}
-\nabla\cdot\left(\overline{A}\left(x\right)\nabla u\left(x\right)\right)=f\left(x\right)\ \text{in}\ \Omega.\label{eq:intro-homogenized-pde}
\end{equation}
\sepia{The homogenized tensor is given by $\overline{A}\left(x\right)=\int_{\bT^{d}}\left(A\left(x,\lambda\right)+A\left(x,\lambda\right)\nabla\chi\left(\lambda;x\right)^{T}\right)\d{\lambda}$, with $\bT^{d}:=\bR^{d}/\bZ^{d}$ being the $d$-dimensional torus, and $\chi\left(\cdot;x\right)$ is known} as the solution the cell problem
\begin{equation}
-\nabla_{\lambda}\cdot\left(\nabla_{\lambda}\chi A^{T}\right)=\nabla_{\lambda}\cdot A^{T},\chi\left(\cdot;x\right)\ \text{is 1-periodic}.\label{eq:intro-cell-problem}
\end{equation}
Solving Eqn. \ref{eq:intro-cell-problem} is expensive even if the underlying space is 1-dimensional, since $\chi\left(\cdot;x\right)$ must be obtained for each $x$. Moreover, this analytical approach may fail in cases where the exact form of the homogenization (Eqn. \ref{eq:intro-homogenized-pde} and \ref{eq:intro-cell-problem}) is out of reach.

Due to the aforementioned reasons, numerical homogenization is often preferred. The heterogeneous multiscale method (HMM) \cite{Abdulle2012,E2003} is a general methodology capable of handling variational \cite{Abdulle2004,Yue2007,Abdulle2008,Abdulle2011} and dynamical \cite{Abdulle2003,Li2005,Engquist2005} problems where the coupling between different scales is deployed to provide necessary information for the macroscopic scale. There have been several works addressing the error estimation \cite{WEINAN2004} as well as extensions with use of variable step sizes \cite{Lee2014,Lee2015}. Another approach is based on the variational connection between scales. The variational multiscale method \cite{Hughes1998}, studied by Hughes et al., decomposes the trial and test spaces into coarse and fine components and proposes a subgrid-scale FEM approach with a posterior error estimate. Other feasible methods include, but are not limited to, the equation-free approach \cite{Analysis2003}, projection-based homogenization \cite{Engquist2008,Nolen2008}, and reiterated-specific homogenization \cite{Engquist2013}. It is worth mentioning that offline-online strategies can also be advantageous if the solver is evaluated for various different source terms. The general idea is to fully resolve the finest scale in the offline phase to extract effective information, so one operates on the coarser level only during the online phase. As one of the successful theories, the multiscale finite element method (MsFEM), proposed and developed by \cite{Hou1999,T.Y.Hou1997,Efendiev2009}, aims to capture local oscillations by replacing the canonical basis with a set of multiscale basis which solves a set of local equations with no source term and matching boundary conditions. To overcome scale resonance, one can also adopt the oversampling technique to reduce the influence of boundary effects.

It should be noted that similar techniques have been developed to perform numerical averaging for dynamical systems. FLow AVeraging integratORS (FLAVORS), proposed and studied in \cite{Tao2011a,Tao2010}, is a class of integrators for stiff ODEs and PDEs. FLAVORS is versatile and non-intrusive since it bypasses the need for identifying variables on different scales. The Seamless multiscale method \cite{E2009}, proposed by E, Ren, and Vanden-Eijnden, has the distinct feature that no re-initialization is required for the microscopic scale at each macroscopic time step, thus avoiding the need for time scale separation. A key component of this seamless method is to perform time rescaling on the microscopic equation. The validity of this technique is based on the conversion relationship between the time variable $t$ and the scale variable $\epsilon$. Albeit elementary, the concept of scale manipulation has been widely applied in many different multiscale settings. The well-established Car-Parrinello technique \cite{carUnifiedApproachMolecular1985} for calculating ab-initio molecular dynamics features a larger fictitious electron mass that characterizes the finest scale. In the simulation of turbulent fluid flows, for example the large eddy simulation method \cite{bardinaImprovedTurbulenceModels1983} and the $K$-$\epsilon$ method \cite{jonesPredictionLaminarizationTwoequation1972}, different modes introduce an eddy viscosity that effectively lowers the Reynolds number, thus enabling a full simulation of homogeneous turbulent flows. As a similar example, the artificial compressibility (AC) \cite{chorinNumericalMethodSolving1967} allows slight perturbation to the incompressibility condition within a certain tolerance, making it possible to apply methods for compressible flows and enabling fast computation; as an analogy, the AC parameter indicates on which scale the compressibility condition starts to be violated.

These methods share the common feature that the rescaling operation is applied to the key components which is closely connected to the scale parameter. However, in the Seamless and FLAVORS method, this rescaling transform applies to the time variable only and there is little attention on adopting rescaling on elliptic systems. Motivated by this, we propose to manipulate the scale in oscillatory elliptic operators to achieve a similar effect. Benefiting from this, the proposed dilation method has its formal simplicity in that it can be used as a drop-in replacement in the existing system, as long as a similar homogenization theory can be established. Besides, there are no overheads such as computing basis functions or homogenized tensors. 

The outline of the paper is listed as follows. We show how to apply numerical averaging techniques on 2-point boundary value problems in Sec. \ref{sec:Reformulation-as-a-time-averaging-problem}. The key connection between the BVP and the corresponding dynamical system is characterized by Thm. \ref{thm:equivalence between BVP and dynamical system} in Sec. \ref{subsec:Rephrasing-the-second} We discuss how to efficiently solve the IVP by numerical averaging in Sec. \ref{subsec:Time-averaging-techniques}. After that, in Sec. \ref{sec:Dilation-Solver}, we study the proposed dilation solver and its properties. The motivation arises from the view of a time dependent system and empirical mode decomposition (Sec. \ref{subsec:Averaging-and-partial-dilation}). \sepia{The rigorous definition of the dilated system is shown} in Sec. \ref{subsec:Dilation-operator}, followed by a proof on the error estimate (Thm. \ref{cor:fem-solver-full-error}) including discussion on numerical schemes (Sec. \ref{subsec:numerical-aspects-and-concerns}). A few numerical examples and corresponding error analysis are covered in Sec. \ref{subsec:experiments} to show the feasibility and stability of the proposed method. Besides, we also study an extension in Sec. \ref{subsec:Structure-aware-dilation} that preserves important structures of the coefficient field and the solution by selective dilation.

\section{Averaging theory\label{sec:Reformulation-as-a-time-averaging-problem}}

The goal of this section is to explore the possibility of applying time averaging techniques to the multiscale problem (Eqn. \ref{eq:intro-elliptic-pde}). The idea is to view the spatial variable $x$ as a temporal one. A crucial step is to rewrite the dynamical system into the form 
\begin{equation}
\frac{\d{}}{\d x}W^{\epsilon}=G\left(W^{\epsilon},\Psi^{\epsilon}\right),\frac{\d{}}{\d x}\Psi^{\epsilon}=\frac{1}{\epsilon}F\left(W^{\epsilon},\Psi^{\epsilon}\right)\label{eq:dsrep-flavors}
\end{equation}
which is required by various time averaging techniques. This is achieved by introducing additional variables that are helpful in moving the position of the $\epsilon^{-1}$ term. The proof and derivation on the equivalence between these two forms is shown in the following paragraph. \sepia{The theoretical foundation for applying numerical averaging methods to Eqn. \ref{eq:dsrep-flavors} is stated as follows.}
\begin{prop}[{\cite[Theorem 5.5]{Artstein1999}}]
\label{cit:averaging theory}Under mild conditions and given $\epsilon$ small enough, the parametrized family of solutions $\left(W^{\epsilon},\Psi^{\epsilon}\right)$ to Eqn. \ref{eq:dsrep-flavors} exists on $x\in\left[0,1\right]$. Furthermore, with 
\[
\mu^{\epsilon}\left(B;T\right):=\frac{1}{T}m\left(\left\{ \tau:0\le\tau\le T,\Psi\left(\epsilon\tau\right)\in B\right\} \right)
\]
being the Dirac distribution along the trajectory and as $\epsilon\to0$, \sepia{$\left\{ \mu^{\epsilon}\right\} $ admits a limit $\mu$ in the sense of weak convergence} and $\left\{ W^{\epsilon}\right\} $ uniformly converges to a process $W$ satisfying $W'=\int G\left(W,\Psi\right)\mu\left(\d{\Psi}\right)$.
\end{prop}
\

In the 1-D setting, we assume $\Omega=\left[0,1\right]$ without loss of generality. We focus on the initial value problem with an assumed condition on the derivative: 
\begin{equation}
-\frac{\d{}}{\d x}\left(A^{\epsilon}\left(x\right)\frac{\d{u^{\epsilon}}}{\d x}\right)=f\left(x,u^{\epsilon}\right),u^{\epsilon}\left(0\right)=0,\frac{\d{}}{\d x}u^{\epsilon}\left(0\right)=v_{0}\label{eq:1d-elliptic-full-equation}
\end{equation}
where \sepia{the real-valued diffusion coefficient takes the form} $A^{\epsilon}\left(x\right)=A\left(x,\frac{x}{\epsilon}\right)$ and $A$ is 1-periodic in the second argument. To ensure well-posedness, we assume $A$ is positive and uniformly bounded from below. The boundary value problem is then addressed by applying the shooting method to Eqn. \ref{eq:1d-elliptic-full-equation}.

\subsection{Rephrasing the second order equation\label{subsec:Rephrasing-the-second}}

A general difficulty lies in the fact that the $\epsilon$ scale occurs in the variable coefficient $A\left(\cdot,\frac{x}{\epsilon}\right)$ rather than in the coefficient as in $\frac{1}{\epsilon}F\left(W,\Psi\right)$. As an illustration, consider a flux function
\begin{equation}
v^{\epsilon}\left(x\right):=\left(2+\sin\left(\frac{2\pi x}{\epsilon}\right)\right)\dfrac{\d{}}{\d x}u^{\epsilon}\left(x\right).\label{eq:equivalence-example-vu}
\end{equation}
We wish to obtain a set of expressions where $\frac{1}{\epsilon}$ appears as an leading coefficient. Even though there is no way to directly extract $\frac{1}{\epsilon}$ from within the $\sin$ function, this can be achieved by introducing more variables. Let $z^{\epsilon}\left(x\right):=\sin\left(\frac{2\pi x}{\epsilon}\right)$; then the derivative of $z^{\epsilon}$ reads $\frac{2\pi}{\epsilon}\cos\left(\frac{2\pi x}{\epsilon}\right)$. To get rid of the remaining $\cos$ term, let $y^{\epsilon}\left(x\right):=\cos\left(\frac{2\pi x}{\epsilon}\right)$ and then the derivative of $\left(z^{\epsilon},y^{\epsilon}\right)$ follows 
\[
\dfrac{\d{}}{\d x}z^{\epsilon}=\frac{2\pi}{\epsilon}y^{\epsilon},\dfrac{\d{}}{\d x}y^{\epsilon}=-\frac{2\pi}{\epsilon}z^{\epsilon}.
\]
Meanwhile, Eqn. \ref{eq:equivalence-example-vu} now reads 
\[
\dfrac{\d{}}{\d x}u^{\epsilon}\left(x\right)=\frac{v^{\epsilon}\left(x\right)}{2+z^{\epsilon}\left(x\right)}
\]
and the $\frac{1}{\epsilon}$ term no longer appears in the frequency domain. We apply this technique in the general setting of the second order equation (Eqn. \ref{eq:1d-elliptic-full-equation}).
\begin{thm}
\label{thm:equivalence between BVP and dynamical system}There exist two multi-variate $C^{1}$-functions $F$ and $G$, independent of $\epsilon$, such that for any $\epsilon>0$ and $u^{\epsilon}\in C^{2}\left(\left[0,1\right]\right)$ satisfying 
\begin{equation}
-\frac{\d{}}{\d x}\left(A^{\epsilon}\left(x\right)\frac{\d{u^{\epsilon}}}{\d x}\right)=f\left(x,u^{\epsilon}\right),\label{eq:equivalence-BVP}
\end{equation}
there exists vector-valued $C^{1}$-functions $W^{\epsilon}=\left(w_{1}^{\epsilon},\dots,w_{k}^{\epsilon}\right)^{T},\Psi^{\epsilon}=\left(\psi_{1}^{\epsilon},\dots\psi_{k}^{\epsilon}\right)^{T}$ s.t. 
\begin{equation}
\frac{\d{}}{\d x}W^{\epsilon}=G\left(W^{\epsilon},\Psi^{\epsilon}\right),\frac{\d{}}{\d x}\Psi^{\epsilon}=\frac{1}{\epsilon}F\left(\Psi^{\epsilon}\right)\label{eq:equivalence-dynamical system}
\end{equation}
for $x\in\left[0,1\right]$ with $w_{1}^{\epsilon}=u^{\epsilon}$. On the contrary, given $C^{1}$-functions $W^{\epsilon}$ and $\Psi^{\epsilon}$ satisfying Eqn. \ref{eq:equivalence-dynamical system} under a suitable initial condition, then the first component $w_{1}^{\epsilon}$ is $C^{2}$ and it solves Eqn. \ref{eq:equivalence-BVP}.
\end{thm}
\begin{proof}
The central idea of this proof is to show a possible form of $F$ and $G$ and the correspondence between $u^{\epsilon}$ and $W^{\epsilon}$. Once the bijection is established, the equivalence between the two equations (Eqn. \ref{eq:equivalence-BVP} and \ref{eq:equivalence-dynamical system}) naturally follows.

\textbf{Step 1:} we construct a pair of $F$ and $G$. Let $v^{\epsilon}$ denote $A^{\epsilon}\left(x\right)\frac{\d{u^{\epsilon}}}{\d x}$ and $v^{\epsilon}\in C^{1}$ as a result of the regularity of $u^{\epsilon}$. The dynamics is given by $\frac{\d{}}{\d x}v^{\epsilon}=-f\left(x,u^{\epsilon}\right)$. As motivated in the aforementioned paragraph, the fast variable $x/\epsilon$ lives in the one dimensional torus $\bT^{1}=\bR\backslash\bZ$ which can be embedded to $\bR^{2}$ by introducing $\phi:\lambda\mapsto\left(\cos2\pi\lambda,\sin2\pi\lambda\right)$. Notice that $\phi$ is parametrized by its arc-length, i.e. $\phi:\left[0,1\right]\to\cS^{1}:=\text{Im}\phi$ satisfying $\phi\left(0\right)=\phi\left(1\right)$ and $\abs{\phi'\left(x\right)}=2\pi$ being a non-vanishing constant. Let $y^{\epsilon}:=\cos2\pi x/\epsilon$ and $z^{\epsilon}:=\sin2\pi x/\epsilon$ denote the first and second component of $\phi\left(x/\epsilon\right)$.

The system $\left(y^{\epsilon}\left(x\right),z^{\epsilon}\left(x\right)\right)$ is autonomous since $\frac{\d{}}{\d x}y^{\epsilon}=-z^{\epsilon}/\epsilon$ and $\frac{\d{}}{\d x}z^{\epsilon}=y^{\epsilon}/\epsilon$. Meanwhile, $\widetilde{A}\left(y,z;x\right):=A\left(x,\phi^{-1}\left(y,z\right)\right)$ is a $C^{1}$-function on $\cS^{1}$ due to the 1-periodicity of $A\left(x,x/\epsilon\right)$. To continuously extend $\widetilde{A}$ to $\bR^{2}$, for each $\left(y,z\right)$ not at the origin, we find the intersection $\left(\widehat{y},\widehat{z}\right)$ of the ray crossing $\left(y,z\right)$ and $\cS^{1}$ such that $\left(y,z\right)=r\left(\widehat{y},\widehat{z}\right)$ and define $\widetilde{A}\left(y,z;x\right):=A\left(x,\phi^{-1}\left(\widehat{y},\widehat{z}\right)\right)\frac{2r^{2}}{1+r^{2}}$ while $\widetilde{A}\left(0,0;x\right):=0$ is a removable singularity.

Therefore, if $W^{\epsilon}$ is identified with $\left(u^{\epsilon},v^{\epsilon},x\right)^{T}$ and $\Psi^{\epsilon}$ with $\left(y^{\epsilon},z^{\epsilon}\right)^{T}$, we naturally define that $F\left(\Psi\right):=\left(-\Psi_{2},\Psi_{1}\right)$ and
\[
G\left(W,\Psi\right):=\left(w_{2}/\widetilde{A}\left(\Psi;w_{3}\right),-f\left(w_{3},w_{1}\right),1\right)^{T}.
\]

\textbf{Step 2:} we show that if $u^{\epsilon}$ solves Eqn. \ref{eq:equivalence-BVP}, then $W^{\epsilon}:=\left(u^{\epsilon},v^{\epsilon},x\right)$ and $\Psi^{\epsilon}:=\left(y^{\epsilon},z^{\epsilon}\right)$ solves Eqn. \ref{eq:equivalence-dynamical system} where $v^{\epsilon}:=A^{\epsilon}\left(x\right)\frac{\d{u^{\epsilon}}}{\d x},y^{\epsilon}:=\phi_{1}\left(x/\epsilon\right),z^{\epsilon}:=\phi_{2}\left(x/\epsilon\right)$. In fact, 
\[
\frac{\d{w_{1}^{\epsilon}}}{\d x}=\frac{\d{u^{\epsilon}}}{\d x}=\frac{1}{A^{\epsilon}\left(x\right)}v^{\epsilon}=\frac{1}{A\left(x,x/\epsilon\right)}w_{2}^{\epsilon}=\frac{1}{\widetilde{A}\left(y^{\epsilon},z^{\epsilon};x\right)}w_{2}^{\epsilon}=\frac{1}{\widetilde{A}\left(\psi_{1}^{\epsilon},\psi_{2}^{\epsilon};w_{3}^{\epsilon}\right)}w_{2}^{\epsilon},
\]
\[
\frac{\d{w_{2}^{\epsilon}}}{\d x}=\frac{\d{}}{\d x}\left[A^{\epsilon}\left(x\right)\frac{\d{u^{\epsilon}}}{\d x}\right]=-f\left(x,u^{\epsilon}\right)=-f\left(w_{3}^{\epsilon},w_{1}^{\epsilon}\right),
\]
\[
\frac{\d{\psi_{1}^{\epsilon}}}{\d x}=\frac{\d{y^{\epsilon}}}{\d x}=\frac{1}{\epsilon}\partial_{1}\phi\left(x/\epsilon\right)=\frac{1}{\epsilon}\partial_{1}\phi\left(\phi^{-1}\left(\psi_{1}^{\epsilon},\psi_{2}^{\epsilon}\right)\right).
\]
Similarly, $\frac{\d{\psi_{2}^{\epsilon}}}{\d x}=\frac{1}{\epsilon}\partial_{2}\phi\left(\phi^{-1}\left(\psi_{1}^{\epsilon},\psi_{2}^{\epsilon}\right)\right)$ and $\frac{\d{w_{3}^{\epsilon}}}{\d x}=1$.

\textbf{Step 3:} conversely, we show that if $W^{\epsilon}$ and $\Psi^{\epsilon}$ solves Eqn. \ref{eq:equivalence-dynamical system} with $w_{3}^{\epsilon}\left(0\right)=0,\psi_{1}^{\epsilon}\left(0\right)=\phi_{1}\left(0\right),\psi_{2}^{\epsilon}\left(0\right)=\phi_{2}\left(0\right)$, then $u^{\epsilon}:=w_{1}^{\epsilon}$ is $C^{2}$ and solves Eqn. \ref{eq:equivalence-BVP}. In fact, $\Psi^{\epsilon}\left(x\right)=\phi\left(x/\epsilon\right)$ solves the initial value problem $\frac{\d{}}{\d x}\Psi^{\epsilon}=\frac{1}{\epsilon}\nabla\phi\left(\phi^{-1}\left(\Psi^{\epsilon}\right)\right),\Psi^{\epsilon}\left(0\right)=\phi\left(0\right)$, leading to $\widetilde{A}\left(\Psi^{\epsilon};w_{3}^{\epsilon}\right)=A\left(x,x/\epsilon\right)=A^{\epsilon}\left(x\right)$. Since $w_{2}^{\epsilon}$ is $C^{1}$ and $\frac{\d{}}{\d x}w_{1}^{\epsilon}=w_{2}^{\epsilon}/A^{\epsilon}\left(x\right)$, thus $w_{1}^{\epsilon}$ is $C^{2}$ and $-\frac{\d{}}{\d x}\left[A^{\epsilon}\left(x\right)\frac{\d{w_{1}^{\epsilon}}}{\d x}\right]=-\frac{\d{w_{2}^{\epsilon}}}{\d x}=f\left(x,w_{1}^{\epsilon}\right)$.
\end{proof}
\begin{rem*}
\sepia{The choice of $\phi$ in the proof above can be generalized to any parameterized Jordan curve with bounded derivatives, provided that the dynamics of the two components of $\phi$ have an autonomous closed form. In this case, we utilize the trigonometric functions $\cos$ and $\sin$ to remain consistent with the derivation presented in Eqn. \ref{eq:equivalence-example-vu}.}
\end{rem*}
As a consequence of Thm. \ref{cit:averaging theory} and \ref{thm:equivalence between BVP and dynamical system}, we derive the averaged equation for Eqn. \ref{eq:1d-elliptic-full-equation} as follows.
\begin{cor}
Let $u^{\epsilon}$ denote the solution to Eqn. \ref{eq:1d-elliptic-full-equation} for each $\epsilon>0$. Then, there exists a subsequence of $u^{\epsilon_{i}}$ and a $C^{2}$-process $u_{0}$ that solves
\[
-\frac{\d{}}{\d x}\left(\overline{A}\left(x\right)\frac{\d u}{\d x}\right)=f\left(x,u\right),u\left(0\right)=0,\frac{\d{}}{\d x}u\left(0\right)=v_{0}
\]
s.t. $u^{\epsilon_{i}}$ and $\frac{\d{}}{\d x}u^{\epsilon_{i}}$ converge uniformly to $u_{0}$ and $\frac{\d{}}{\d x}u_{0}$ on $\left[0,1\right]$ respectively, where $\overline{A}\left(x\right)=\left[\int_{0}^{1}A\left(x,\lambda\right)^{-1}\d{\lambda}\right]^{-1}$ is the homogenized coefficient.
\end{cor}
\begin{proof}
\sepia{It suffices to obtain the effective equations} for $\frac{\d{}}{\d x}v^{\epsilon}=-f\left(x,u^{\epsilon}\right)$ and $\frac{\d{}}{\d x}u^{\epsilon}=v^{\epsilon}/\widetilde{A}\left(y^{\epsilon},z^{\epsilon};x\right)$. In fact, the invariant distribution $\mu$ is the Dirac measure along the curve $\phi\left(\left[0,1\right]\right)$, weighted by the arc-length up to a normalization constant. \sepia{Thus, effective source terms are given by $\int-f\left(x,u\right)\d{\mu\left(y,z\right)}=-f\left(x,u\right)$ and
\[
\int_{\phi\left(\left[0,1\right]\right)}\frac{v\left(x\right)\d{\mu\left(y,z\right)}}{\widetilde{A}\left(y,z;x\right)}=\int_{0}^{1}\frac{v\left(x\right)}{\widetilde{A}\left(\phi_{1}\left(\lambda\right),\phi_{2}\left(\lambda\right);x\right)}\d{\lambda}=\int_{0}^{1}\frac{v\left(x\right)}{A\left(x,\lambda\right)}\d{\lambda}=\frac{v\left(x\right)}{\overline{A}\left(x\right)},
\]
leading to the effective equations $\frac{\d{}}{\d x}v=-f\left(x,u\right),\frac{\d{}}{\d x}u=v\left(x\right)/\overline{A}\left(x\right)$.}
\end{proof}

\subsection{Numerical averaging techniques\label{subsec:Time-averaging-techniques}}

We briefly discuss a few time averaging techniques which can be useful for efficient evaluation of Eqn. \ref{eq:equivalence-dynamical system} (equivalently for Eqn. \ref{eq:1d-elliptic-full-equation}), such as Seamless method, FLAVORS, etc. The error analysis is postponed to a future paper with a systematic comparison.

\paragraph*{Seamless method\label{subsec:Seamless-method}}

In \cite{E2009}, a general strategy is proposed for designing solvers for multi-scale problems that aims to effectively increase the $\epsilon$ parameter. With Thm. \ref{thm:equivalence between BVP and dynamical system} kept in mind, we apply the forward Euler scheme on Eqn. \ref{eq:equivalence-dynamical system}
\begin{align}
W_{n+1} & =W_{n}+G\left(W_{n},\Psi_{n}\right)\Delta x,\label{eq:motivation-discretization-macro}\\
\Psi_{n+1} & =\Psi_{n}+\frac{1}{\epsilon}F\left(\Psi_{n}\right)\Delta x.\label{eq:motivation-discretization-micro}
\end{align}
If one uses different clocks for the two scales, for example $\Delta x$ for macroscale (Eqn. \ref{eq:motivation-discretization-macro}) and $\tau<\Delta x$ for microscale (Eqn. \ref{eq:motivation-discretization-micro}), then the modified scheme
\begin{align*}
\widetilde{W}_{n+1} & =\widetilde{W}_{n}+G\left(\widetilde{W}_{n},\widetilde{\Psi}_{n}\right)\Delta x,\\
\widetilde{\Psi}_{n+1} & =\widetilde{\Psi}_{n}+\frac{1}{\epsilon}F\left(\widetilde{\Psi}_{n}\right)\tau=\widetilde{\Psi}_{n}+\frac{\tau}{\epsilon\Delta x}F\left(\widetilde{\Psi}_{n}\right)\Delta x
\end{align*}
can be viewed as the forward Euler scheme applied on the system
\begin{equation}
\frac{\d{}}{\d x}\widetilde{W}=G\left(\widetilde{W},\widetilde{\Psi}\right),\frac{\d{}}{\d x}\widetilde{\Psi}=\frac{\tau}{\epsilon\Delta x}F\left(\widetilde{\Psi}\right)\label{eq:seamless-dynamical-system-dilated}
\end{equation}
where the effective parameter satisfies $\frac{1}{\widetilde{\epsilon}}=\frac{\tau}{\epsilon\Delta x}$, namely $\widetilde{\epsilon}=\frac{\Delta x}{\tau}\epsilon>\epsilon$. Thus, using a faster clock on the fast dynamics effectively increases the $\epsilon$ parameter and reduces the stiffness of the system. We point out that the second order system corresponding to Eqn. \ref{eq:seamless-dynamical-system-dilated} shares the same form as Eqn. \ref{eq:1d-elliptic-full-equation} except for the diffusion coefficient
\begin{equation}
A_{\text{Seamless}}^{\epsilon}\left(x\right)=A\left(x,\frac{\tau}{\epsilon\Delta x}x\right).\label{eq:seamless effective}
\end{equation}
A visualization is displayed in Fig. \ref{fig:seamless sub a}.

\paragraph*{FLAVORS}

First introduced in \cite{Tao2010,Tao2011a}, FLAVORS does not require an explicit identification of the slow and fast variables. Let $\widetilde{W}$ denote the augmented collection of variables $\widetilde{W}:=\left(u,v,x,y,z\right)$. Let $\Phi_{x}^{1/\epsilon}$ denote the solution map $\widetilde{W}^{\epsilon}\left(0\right)\mapsto\widetilde{W}^{\epsilon}\left(x\right)$ under Eqn. \ref{eq:equivalence-dynamical system}. According to FLAVORS, the numerical solution map reads $\Psi_{k\Delta x}:=\left(\Phi_{\Delta x-\tau}^{0}\circ\Phi_{\tau}^{1/\epsilon}\right)^{k}$ where the fast dynamics switches between on ($\Phi_{\tau}^{1/\epsilon}$) and off ($\Phi_{\Delta x-\tau}^{0}$) in an alternating manner so that the fast dynamics is sampled evenly with a relatively lower cost due to $\tau\ll\Delta x$. 

\subsection{Connection between 1-D boundary and initial value problems}

Thm. \ref{thm:equivalence between BVP and dynamical system} provides a possibility of translating a second order equation into a first order one with coefficients independent of $\epsilon$; furthermore, it provides the correspondence between $u^{\epsilon}$ and $\left(W^{\epsilon},\Psi^{\epsilon}\right)$ via $W^{\epsilon}=\left(u^{\epsilon},v^{\epsilon},x\right)^{T},\Psi^{\epsilon}=\left(y^{\epsilon},z^{\epsilon}\right)^{T}$. Since the initial value problem has been addressed, we can now turn to the boundary value problem. The general idea is to try different initial values $v_{0}$ and pick the best one based on how close the terminal point $u\left(1\right)$ of each trajectory is. To efficiently solve the algebraic equation for $v_{0}$, the Newton's method is preferred at the cost of an additional dynamical system that characterizes the partial derivative of the solution map to the assumed initial value $v_{0}$.

\section{Dilation Solver\label{sec:Dilation-Solver}}

In the previous section, we demonstrate how to apply numerical averaging techniques to multiscale boundary value problems. \sepia{Although the translated dynamical system depends on a one-dimensional assumption, the rescaling ideas and techniques can be generalized to a broader setting, allowing for comparison of three possible approaches. The first approach leverages explicit knowledge of the decomposition into faster and slower scales, enabling the direct derivation of a practical algorithm with straightforward error analysis. The second approach utilizes existing pre-processing techniques to detect fast and slow frequencies, which are then used to apply the first technique. The third approach does not require prior knowledge of the scale decomposition; instead, the entire variable coefficient is locally dilated.}

The roadmap is laid as follows. We quickly review the homogenization theory and quantitative estimate in \ref{subsec:Quantitative-stimate-on-homogenization-error}. We then derive the method of partial dilation and explain its connection to averaging techniques in Sec. \ref{subsec:Averaging-and-partial-dilation}, corresponding to the second approach mentioned above. To address the third approach, we propose a local dilation method in Sec. \ref{subsec:Dilation-operator} to circumvent the need for scale identification. The numerical implementation is discussed in Sec. \ref{subsec:numerical-aspects-and-concerns} and the discretization error is analyzed. Finally, we show the numerical performance on two typical examples in Sec. \ref{subsec:experiments}.

\subsection{Quantitative estimate on homogenization error\label{subsec:Quantitative-stimate-on-homogenization-error}}

\sepia{As our focus has shifted to a general dimension, the notation $\Omega$ will now refer to a given bounded open set in $\bR^{d}$.} The following conditions are often assumed to satisfy the Lax-Milgram theorem.
\begin{assumption}
\label{assu:conditions on A}We assume the following three conditions:
\begin{enumerate}
\item \sepia{$A^{\epsilon}\left(x\right)\in M^{d\times d}$ is symmetric for any $x\in\Omega$, where $M^{d\times d}$ denotes the set of $d$-by-$d$ real matrices};
\item $A^{\epsilon}\left(x\right)$ is uniformly elliptic in $\Omega$, i.e. there exists $r\in\left(0,1\right)$ s.t. $\left\langle A^{\epsilon}\left(x\right)\xi,\xi\right\rangle \ge r\abs{\xi}^{2}$ holds for any $x\in\Omega,\xi\in\bR^{d}$;
\item $A^{\epsilon}\left(x\right)$ is bounded, i.e. $\abs{\left\langle A^{\epsilon}\left(x\right)\xi,\eta\right\rangle }\le M\abs{\xi}\,\abs{\eta}$ holds for any $x\in\Omega,\xi,\eta\in\bR^{d}$.
\end{enumerate}
\end{assumption}
When the oscillation only involves the fast variable $x/\epsilon$, a quantitative approximation rate in the $L^{2}\left(\Omega\right)$ space can be concluded from the following result.
\begin{prop}[{\cite[Theorem 1.1]{Xu2016}}]
\label{thm:homogenization-error-Lp-estimate}Suppose that $A^{\epsilon}\left(x\right)=A\left(x/\epsilon\right)$ satisfies Assump. \ref{assu:conditions on A}. Let $u^{\epsilon}$ and $u_{0}$ denote solution to Eqn. \ref{eq:intro-elliptic-pde} and \ref{eq:intro-homogenized-pde}. Then there exists constants $C\left(\Omega,r,M,d\right)$ and $r_{0}\left(\Omega,r,M,d\right)$ s.t. $\norm{u^{\epsilon}-u_{0}}_{L^{2}\left(\Omega\right)}\le C\epsilon\ln\left(r_{0}/\epsilon\right)\norm f_{L^{2}\left(\Omega\right)}$. Furthermore, if $u_{0}\in H^{2}\left(\Omega\right)$, then $\norm{u^{\epsilon}-u_{0}}_{L^{p}\left(\Omega\right)}\le C\epsilon\norm{u_{0}}_{H^{2}\left(\Omega\right)}$ where $p:=\frac{2d}{d-1}$.
\end{prop}
\begin{rem*}
We apply Thm. \ref{thm:homogenization-error-Lp-estimate} when the diffusion tensor depends on both $x$ and $x/\epsilon$ in the later paragraph. This can be partially justified by the fact that the homogenization effect for linear operators is local. We refer to the ``shifted method'' \cite{Zhikov2016} for an alternative approach to prove the quantitative rate of convergence.
\end{rem*}

\subsection{Averaging and partial dilation\label{subsec:Averaging-and-partial-dilation}}

Motivated by the Seamless method in the regime of dynamical systems, we seek to find a similar method that works smoothly for boundary value problems, especially in higher dimensions. In Sec. \ref{subsec:Seamless-method}, we show the equivalence of altering the clock for the microscale and changing the effective scale parameter. To put differently, the Seamless method relaxes a system in the form of $\dot{u}=\cL^{\epsilon}u$ into $\cL^{\widetilde{\epsilon}}$ in the sense that $\widetilde{\epsilon}>\epsilon$. This encourages a direct relaxation of the scale parameter in Eqn. \ref{eq:intro-elliptic-pde}, i.e. to change $A^{\epsilon}$ into $A^{\widetilde{\epsilon}}$.

To put this rigorously, we layout the following preparations. We consider the family of matrix-valued functions that are 1-periodic in the second argument:
\[
\cS:=\left\{ A\in L^{\infty}\left(\Omega\times\bR^{d};M^{d\times d}\right);A\left(x,\lambda+e_{i}\right)=A\left(x,\lambda\right),\forall1\le i\le d\right\} 
\]
where $e_{i}=\left(0,\dots,\underset{i\text{-th}}{1},\dots,0\right)$ is the set of standard basis. Then, oscillatory coefficients can be obtained from elements in $\cS$ via an $\epsilon$-wrapping operation, \sepia{defined as
\begin{align}
\iota^{\epsilon}:\cS & \to L^{\infty}\left(\Omega;M^{d\times d}\right),\label{eq:partial dilation-wrapping def}\\
A & \mapsto A^{\epsilon},A^{\epsilon}\left(x\right):=A\left(x,x/\epsilon\right)\nonumber 
\end{align}
}and let $\cT^{\epsilon}:=\text{Im}\iota^{\epsilon}$ denote all possible oscillatory coefficients. The partial dilation operator is defined as
\begin{equation}
\cD_{m}^{\text{p}}\left(B\right)\left(x,\lambda\right):=B\left(x,\frac{\lambda}{m}\right),B\in\cS.\label{eq:partial dilation-partial dilation op def}
\end{equation}

We wish to derive a relaxed coefficient $A^{m\epsilon}=\iota^{m\epsilon}A$ from a rather oscillatory one $A^{\epsilon}=\iota^{\epsilon}A$. In fact, by definitions in Eqn. \ref{eq:partial dilation-wrapping def} and \ref{eq:partial dilation-partial dilation op def}, we have
\[
\left(\iota^{m\epsilon}A\right)\left(x\right)=A\left(x,x/m\epsilon\right)=\left(\cD_{m}^{\text{p}}A\right)\left(x,x/\epsilon\right)=\left(\iota^{\epsilon}\cD_{m}^{\text{p}}A\right)\left(x\right),
\]
thus the partially dilated coefficient can be represented as
\[
A^{m\epsilon}=\iota^{m\epsilon}A=\iota^{\epsilon}\cD_{m}^{\text{p}}A=\iota^{\epsilon}\cD_{m}^{\text{p}}\left(\iota^{\epsilon}\right)^{-1}A^{\epsilon}.
\]
We thus define the partially dilated problem to Eqn. \ref{eq:intro-elliptic-pde} as
\begin{align}
-\nabla\cdot\left[\iota^{\epsilon}\cD_{m}^{\text{p}}\left(\iota^{\epsilon}\right)^{-1}A^{\epsilon}\left(x\right)\nabla u^{\epsilon}\left(x\right)\right] & =f\left(x\right) & \text{in}\ \Omega,\label{eq:partial dilation-dilated-pde}\\
u^{\epsilon} & =0 & \text{on}\ \partial\Omega\nonumber 
\end{align}
where Eqn. \ref{eq:partial dilation-dilated-pde} can also be written as $-\nabla\cdot\left[A\left(x,x/m\epsilon\right)\nabla u^{\epsilon}\left(x\right)\right]=f\left(x\right)$ if the scale identification $A^{\epsilon}\left(x\right)=A\left(x,x/\epsilon\right)$ is known. We include the following diagram to illustrate the connection between $A^{m\epsilon},A^{\epsilon}$, and $A$.

\[\begin{tikzcd}[ampersand replacement=\&]
A^{\epsilon}
\arrow[r, "\left(\iota^{\epsilon}\right)^{-1}"]
\arrow[d, "\iota^{\epsilon}\cD_{m}^{\text{p}}\left(\iota^{\epsilon}\right)^{-1}"']
\&
A
\arrow[r, "\cH"]
\&
\overline {A}
\\
A^{m\epsilon}
\arrow[ru, "\left(\iota^{m\epsilon}\right)^{-1}"']
\&
\&
\end{tikzcd}\]

\subsubsection{Error analysis\protect\footnote{Notation wise, the gradient of a vector $v\left(z\right)$ is $\left(\nabla v\right)_{ij}:=\frac{\partial v_{i}}{\partial z_{j}}$ and the divergence of a matrix $A\left(z\right)$ reads $\left(\nabla\cdot A\right)_{i}:=\sum_{j}\frac{\partial A_{ij}}{\partial z_{j}}$. The Euclidean norm $\protect\abs z$ is the 2-norm for a vector $z$ in $\protect\bR^{d}$. By default, the symbol $x$ refers to a vector in the macroscopic domain $\Omega\subset\protect\bR^{d}$ while $\lambda$ refers to a vector in the microscopic domain $\protect\bT^{d}$ ($d$-dimensional torus).}}

To derive the error estimate on Eqn. \ref{eq:partial dilation-dilated-pde}, the homogenized equation is introduced to bridge the gap. Let the homogenized operator be
\begin{equation}
\cH:\cS\to L^{\infty}\left(\Omega;M^{d\times d}\right),A\mapsto\overline{A}:=\int_{\bT^{d}}\left(A\left(x,\lambda\right)+A\left(x,\lambda\right)\nabla_{\lambda}\chi\left(\lambda;x\right)^{T}\right)\d{\lambda}\label{eq:homogenization-operator-def}
\end{equation}
where $\chi:\Omega\to\left[C_{per}^{1}\left(\bT^{d}\right)\right]^{d}$ solves a family of cell problems (Eqn. \ref{eq:intro-cell-problem}). 

Let $u_{0}$ denote the solution to the homogenized equation (Eqn. \ref{eq:intro-homogenized-pde}), $u_{\cD^{\text{p}}}^{\epsilon}$ to Eqn. \ref{eq:partial dilation-dilated-pde}, and $u_{\cD^{\text{p}},h}^{\epsilon}\in V_{h}$ the piecewise-linear FEM solution to Eqn. \ref{eq:partial dilation-dilated-pde} (more details on FEM are in Sec. \ref{subsec:numerical-aspects-and-concerns}). By the triangle inequality, 
\[
\norm{u_{0}-u_{\cD^{\text{p}},h}^{\epsilon}}\le\underbrace{\norm{u_{0}-u_{\cD^{\text{p}}}^{\epsilon}}}_{\text{homogenization error}}+\underbrace{\norm{u_{\cD^{\text{p}}}^{\epsilon}-u_{\cD^{\text{p}},h}^{\epsilon}}}_{\text{discretization error}}.
\]

\paragraph{Homogenization error}

Since $u_{\cD^{\text{p}}}^{\epsilon}$ solves \ref{eq:partial dilation-dilated-pde} and $\iota^{\epsilon}\cD_{m}^{\text{p}}\left(\iota^{\epsilon}\right)^{-1}A^{\epsilon}\left(x\right)$ turns out to be $A\left(x,x/m\epsilon\right)$, we can directly apply Prop. \ref{thm:homogenization-error-Lp-estimate} to show that
\[
\norm{u_{0}-u_{\cD^{\text{p}}}^{\epsilon}}_{L^{2}\left(\Omega\right)}\le\norm{u_{0}-u_{\cD^{\text{p}}}^{\epsilon}}_{L^{2d/\left(d-1\right)}\left(\Omega\right)}\le C\left(\Omega,r,M,d\right)m\epsilon\norm{u_{0}}_{H^{2}\left(\Omega\right)}.
\]

\paragraph{Discretization error}

We apply the FEM error estimate \cite{Ciarlet1991} to show that
\[
\norm{u_{\cD^{\text{p}}}^{\epsilon}-u_{\cD^{\text{p}},h}^{\epsilon}}_{L^{2}\left(\Omega\right)}\le C\left(\widehat{T},\widehat{P},\widehat{\Sigma}\right)h^{2}\norm{u_{\cD^{\text{p}}}^{\epsilon}}_{H^{2}\left(\Omega\right)}
\]
under Assump. \ref{assu:conditions on A}, where $\left(\widehat{T},\widehat{P},\widehat{\Sigma}\right)$ stands for the triplet of the reference finite element. \sepia{Since $A^{m\epsilon}=\iota^{\epsilon}\cD_{m}^{\text{p}}\left(\iota^{\epsilon}\right)^{-1}A^{\epsilon}$, the effective scale parameter in Eqn. \ref{eq:partial dilation-dilated-pde} is $m\epsilon$. Motivated by the two-scale asymptotic ansatz \cite[Section 12.4]{Grigorios2008} 
\[
u_{\cD^{\text{p}}}^{\epsilon}\left(x\right)=u_{0}\left(x\right)+m\epsilon u_{1}\left(x,x/\left(m\epsilon\right)\right)+m^{2}\epsilon^{2}u_{2}\left(x,x/\left(m\epsilon\right)\right)+\text{h.o.t.},
\]
an informed conjecture is that the $H^{2}$-norm of $u_{\cD^{\text{p}}}^{\epsilon}$ is dominated by $\norm{\epsilon^{-1}\partial_{\lambda}^{2}u_{1}}_{L^{2}}$. This can be made rigorous by matching the power terms of $\epsilon$ and recursively defining equations for $u_{2}$ and other higher-order terms \cite[Section 2.5]{Bensoussan1978}. Consequently, the FEM solution error is bounded by}
\[
\norm{u_{\cD^{\text{p}}}^{\epsilon}-u_{\cD^{\text{p}},h}^{\epsilon}}_{L^{2}\left(\Omega\right)}\le\frac{1}{\epsilon}C\left(\widehat{T},\widehat{P},\widehat{\Sigma}\right)\norm{\partial_{\lambda}^{2}u_{1}}_{L^{2}\left(\Omega\times\bT^{d}\right)}\norm{\chi}_{C\left(\Omega;H^{2}\left(\Pi\right)\right)}h^{2}.
\]
In the case where $A$ and $f$ have better regularity, it is also possible to adopt basis functions with higher orders so that it leads to faster convergence. For example, if $u_{\cD^{\text{p}}}^{\epsilon}\in H^{k+1}\left(\Omega\right),$then the approximation error reads $\cO\left(h^{k+1}\right)\norm{u_{\cD^{\text{p}}}^{\epsilon}}_{H^{k+1}\left(\Omega\right)}$ where $\norm{u_{\cD^{\text{p}}}^{\epsilon}}_{H^{k+1}\left(\Omega\right)}\apprle\epsilon^{-k}\norm{\partial_{\lambda}^{k+1}u_{1}}_{L^{2}\left(\Omega\times\bT^{d}\right)}$, implying that the overall bound is $\cO\left(h^{k+1}\epsilon^{-k}\right)$.

To summarize, the error of the FEM solution to the homogenized solution is given by $\cO\left(m\epsilon+h^{k+1}\epsilon^{-k}\right)$ for some $k\ge1$\@. This is compared with other methods in two numerical experiments in Sec. \ref{subsec:integrated-test}.

\subsubsection{Numerical inversion $\left(\iota^{\epsilon}\right)^{-1}$\label{subsec:Numerical-inversion}}

Compared to the local dilation method in the following paragraph, this partial dilation method is much easier to comprehend and to analyze. However, a major requirement is that the scale separated form $A\left(x,x/\epsilon\right)$ must be known while solving Eqn. \ref{eq:partial dilation-dilated-pde}. Equivalently, we either reduce the inverse operator $\left(\iota^{\epsilon}\right)^{-1}$ or solve it numerically.

Since $\left(\iota^{\epsilon}\right)^{-1}$ maps $A^{\epsilon}\left(x\right)$ back to $A\left(x,\lambda\right),\lambda=x/\epsilon$, it amounts to decompose the oscillations in $A^{\epsilon}\left(x\right)$ into a few periodic mode functions. To this end, signal decomposition techniques such as Empirical Mode Decomposition (EMD \cite{Xu2006}) and Synchrosqueezed Wavelet Transforms (SWT \cite{Daubechies2009}) can be employed to extract oscillatory components and modify their instantaneous frequencies. This enables an extension of the partial dilation method that falls in the second category, where we do not require an explicit scale identification which is inferred by decomposition. We demonstrate the feasibility in the numerical section (Sec. \ref{subsec:integrated-test}) and we intend to further explore the possibility in future works.

\subsection{Dilation operator\label{subsec:Dilation-operator}}

Straightforward as the partial dilation method is, the application is often restricted due to the need of explicit scale separation or the expensive numerical inversion process. It is thus tempting to ask if we can bypass this issue and derive an approximation to $\iota^{\epsilon}\cD_{m}^{\text{p}}\left(\iota^{\epsilon}\right)^{-1}$ directly. To better illustrate the construction of local dilation operators, we begin with an example where $A^{\epsilon}\left(x\right)=A_{0}\left(x\right)+A_{1}\left(x\right)\varphi\left(x/\epsilon\right)$ has a slowing varying component $A_{0}$ and a fast oscillatory component $\varphi\left(x/\epsilon\right)$. Applying the partial dilation operator yields
\[
\left[\iota^{\epsilon}\cD_{m}^{\text{p}}\left(\iota^{\epsilon}\right)^{-1}A^{\epsilon}\right]\left(x\right)=\left(\iota^{\epsilon}\cD_{m}^{\text{p}}\right)\left[A_{0}\left(x\right)+A_{1}\left(x\right)\varphi\left(x/\epsilon\right)\right]=A_{0}\left(x\right)+A_{1}\left(x\right)\varphi\left(x/m\epsilon\right).
\]
Following the idea of rescaling, one might wish to try substituting $x$ by $x/m$, but in general $A_{i}\left(x\right)$ can be quite different from $A_{i}\left(x/m\right)$ on the domain $\Omega$. Thus, we are seeking for a coordinate transform $\phi\left(x\right)$ that achieves $\left[\iota^{\epsilon}\cD_{m}^{\text{p}}\left(\iota^{\epsilon}\right)^{-1}A^{\epsilon}\right]\left(x\right)\approx A^{\epsilon}\left(\phi\left(x\right)\right)$ in the way that
\begin{enumerate}
\item \sepia{every partial derivative $\partial_{x_{i}}\phi$ is close to $1/m$ so that $\varphi\left(\phi\left(x\right)/\epsilon\right)$ and $\varphi\left(x/m\epsilon\right)$, as fast variables, have similar distributions}, while
\item $\phi$ tracks the identity function on the macroscale so that $A_{i}\left(\phi\left(x\right)\right)\approx A_{i}\left(x\right)$.
\end{enumerate}
One such candidate is the family of local dilation operators which relax the instantaneous frequency by shrinking the input coordinates locally. These operators are defined in the following section.

\subsubsection{Shrinkage mapping and dilation operator}
\begin{figure}[tbph]
\begin{centering}
\includegraphics[width=0.9\linewidth]{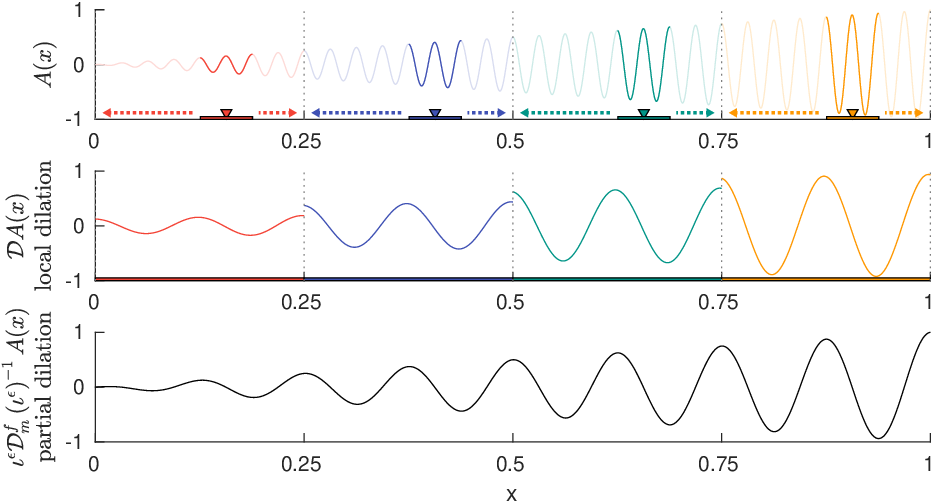}
\par\end{centering}
\caption{A 1-D example comparing the original coefficient $A\left(x\right)=x\cos64\pi x$ and the dilated versions. The parameters of the (local) dilated operator $\protect\cD_{L,m,\nu}$ read $L=\frac{1}{4},m=4,\nu=0.67$. The highlighted area in the top row indicates the range of the shrinkage mapping $\phi$, where each segment is centered at the anchor function $\Phi$. The second row shows the local dilation while the third row is from partial dilation, i.e. the effective coefficient in the Seamless method (Eqn. \ref{eq:seamless effective}).}
\label{fig:seamless sub a}
\end{figure}

\begin{figure}[tbph]
\begin{centering}
\includegraphics[width=0.98\linewidth]{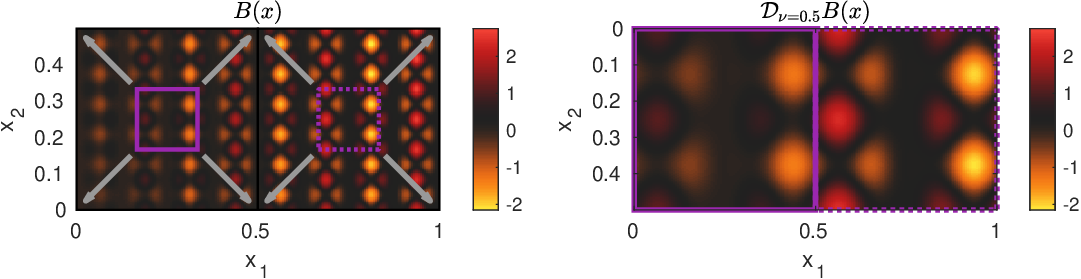}
\par\end{centering}
\caption{A 2-D example of applying dilation operator $\protect\cD_{L=0.5,m=3,\nu=0.5}$ on $B\left(x\right)=\left(2x_{1}-x_{1}^{2}+x_{2}-x_{2}^{2}\right)\left[\sin24\pi x_{1}+\cos24\pi x_{2}\right]\protect\abs{\sin\left(8\pi x_{1}\right)}^{1.5}+0.3x_{1}$. }
\label{fig:seamless sub b}
\end{figure}

\begin{defn}
\label{def:anchor-shrinkage-and-dilation}For a given \sepia{mesoscopic length $L>0$} and an anchoring factor $\nu\in\left[0,1\right]$, the anchor function \sepia{$\Phi_{L,\nu}:\bR\to\bR$} is defined as
\[
\Phi_{L,\nu}\left(y\right):=\left(\left\lfloor \frac{y}{L}\right\rfloor +\nu\right)L
\]
where $\left\lfloor \cdot\right\rfloor $ refers to the floor function. Based on this and a scaling factor $m>1$, the shrinkage mapping is defined as 
\[
\phi_{L,m,\nu}\left(y\right):=\frac{1}{m}\left(y-\Phi_{L,\nu}\left(y\right)\right)+\Phi_{L,\nu}\left(y\right).
\]
Finally, the dilation operator transforms a given function to one that takes a modified input according to
\begin{equation}
\left(\cD_{L,m,\nu}B\right)\left(x_{1},\dots,x_{n}\right):=B\left(\phi_{L,m,\nu}\left(x_{1}\right),\dots,\phi_{L,m,\nu}\left(x_{n}\right)\right).\label{eq:dilation-operator-def}
\end{equation}
\end{defn}
To help understand the operators defined above, we include two examples in Fig. \ref{fig:seamless sub a} and \ref{fig:seamless sub b}.

With a slight abuse of notation, the shrinkage operator $\phi_{L,m,\nu}$ can apply to a vector $x=\left(x_{1},\dots,x_{d}\right)$ in a component-wise manner
\[
\phi_{L,m,\nu}\left(x\right):=\left(\phi_{L,m,\nu}\left(x_{1}\right),\dots,\phi_{L,m,\nu}\left(x_{d}\right)\right).
\]
Thus, Eqn. \ref{eq:dilation-operator-def} has a concise form $\left(\cD_{L,m}B\right)\left(x\right)=B\left(\phi_{L,m}\left(x\right)\right)$. In the following paragraphs, we omit $L$, $m$, and $\nu$ when there is no ambiguity on parameters. We provide an estimate of the distance between $\phi\left(x\right)$ and $x$.
\begin{lem}
\label{lem:shrinkage-difference-estimate}For given $L$, $m$, and $\nu$, $\abs{\phi\left(x\right)-x}\le dL$ holds for any $x\in\Omega$.
\end{lem}
\begin{proof}
We first relax the vector 2-norm by $\abs{\phi\left(x\right)-x}\le\sum_{i=1}^{d}\abs{\phi\left(x_{i}\right)-x_{i}}$. Then, for each index $1\le i\le d$, by definition of the anchor function and shrinkage mapping,
\[
\abs{\phi\left(x_{i}\right)-x_{i}}=\abs{\left(1-\frac{1}{m}\right)\left(x_{i}-\Phi\left(x_{i}\right)\right)}=\frac{m-1}{m}L\abs{\left\{ \frac{x_{i}}{L}\right\} -\nu}.
\]
Since $m>1$, $0\le x_{i}/L<1$, and $0\le\nu\le1$, the RHS above shall be smaller than $L$, leading to $\abs{\phi\left(x\right)-x}\le\sum_{i=1}^{d}\abs{\phi\left(x_{i}\right)-x_{i}}\le dL.$
\end{proof}
\begin{rem*}
\sepia{Lem. \ref{lem:shrinkage-difference-estimate} provides the best estimate of $\phi\left(x\right)-x$ in terms of the exponent of $L$, regardless of the choice of $\nu$. However, the convergence rate of the overall dilation solver might improve with a symmetry-motivated choice of $\nu=\frac{1}{2}$ (see discussion in Par. \ref{par:numerical-dilation-error}).}
\end{rem*}

\subsubsection{Dilated system}

With the aforementioned notations, we propose to solve the following dilated problem
\begin{align}
-\nabla\cdot\left[\cD A^{\epsilon}\left(x\right)\nabla u_{D}^{\epsilon}\left(x\right)\right] & =f\left(x\right) & \text{in}\ \Omega,\label{eq:dilation-solver-dilated-pde}\\
u_{D}^{\epsilon} & =0 & \text{on}\ \partial\Omega
\end{align}
under a proper set of parameters $L$, $m$, and $\nu$. We point out that $\cD A^{\epsilon}\in\cT^{\epsilon}$ if $\Omega$ has a nice geometric structure.
\begin{lem}
\label{lem:DAeps-representation}For given $A\in C^{1}\left(\Omega\times\bT^{d}\right)^{d\times d}$ on $\Omega=\left[0,1\right]^{d}$, let $A^{\epsilon}:=\iota^{\epsilon}A$ denote the wrapped oscillatory coefficient. Then, for a given dilation operator $\cD_{L,m,\nu}$, there exists $\widetilde{A}\in\cS$ s.t. $\cD_{L,m,\nu}A^{\epsilon}=\iota^{m\epsilon}\widetilde{A}$. Furthermore, $\widetilde{A}$ is a $C^{1}\left(T_{k}\times\bT^{d};M^{d\times d}\right)$ matrix-valued function if restricted in $T_{k}$, where $T_{k}:=\prod_{i=1}^{d}\left[\left(k_{i}-1\right)L,k_{i}L\right)$ is a $d$-dimensional cube and $k=\left(k_{1},k_{2},\dots,k_{d}\right)$ is a multi-index with $1\le k_{i}\le\frac{1}{L}$ for $1\le i\le d$.
\end{lem}
\begin{proof}
Let $\widetilde{\epsilon}:=m\epsilon$ denote the relaxed scale parameter. Notice that the family of cubes $\left\{ T_{k}\right\} $ forms a partition of $\Omega$. The difference in the fast variable reads
\[
\frac{\phi\left(x\right)}{\epsilon}-\frac{x}{\widetilde{\epsilon}}=\frac{1}{\epsilon}\left[\frac{1}{m}\left(x-\Phi\left(x\right)\right)+\Phi\left(x\right)\right]-\frac{x}{m\epsilon}=\frac{1}{\epsilon}\left(1-\frac{1}{m}\right)\Phi\left(x\right)
\]
where $\Phi\left(x\right)=\left(\left\lfloor \frac{x}{L}\right\rfloor +\nu\right)L$ is a constant function within each cube $T_{k}$, so the multi-index $k$ is associated to $\varphi_{k}:=\frac{1}{\epsilon}\left(1-\frac{1}{m}\right)\left(k+\nu\right)L$. Then, the dilated coefficient reads
\[
\left(\cD A^{\epsilon}\right)\left(x\right)=A\left(\phi\left(x\right),\frac{\phi\left(x\right)}{\epsilon}\right)=A\left(\phi\left(x\right),\frac{x}{\widetilde{\epsilon}}+\varphi_{k}\right)\ \text{for}\ x\in T_{k}
\]
which naturally leads to $\widetilde{A}\left(x,\lambda\right):=A\left(\phi\left(x\right),\lambda+\varphi_{k}\right)$ for $x\in T_{k}$. $\widetilde{A}$ is 1-periodic in $\lambda$ due to $1$-periodicity of $A$, so $\widetilde{A}\in\cS$ and one can verify that $\cD A^{\epsilon}=\iota^{m\epsilon}\widetilde{A}$. In each $T_{k}$, since $A$ is $C^{1}$ and $\phi|_{T_{k}}$ is affine, $\widetilde{A}|_{T_{k}}$ is $C^{1}$ as well.
\end{proof}
\begin{rem*}
In the context of Lem. \ref{lem:DAeps-representation}, since $\widetilde{A}\in\cS$, we can apply the homogenization operator $\cH$ on it, leading to the homogenized dilated system
\begin{align}
-\nabla\cdot\left[\cH\widetilde{A}\left(x\right)\nabla u_{0,D}\left(x\right)\right] & =f\left(x\right) & \text{in}\ \Omega,\label{eq:dilation-solver-homogenized-dilated-pde}\\
u_{0,D} & =0 & \text{on}\ \partial\Omega.
\end{align}
\end{rem*}
We point out that $\cH$ and $\cD$ are interchangeable in the following sense.

\begin{lem}
\label{lem:interchangability}In the context of Lem. \ref{lem:DAeps-representation}, the dilation of homogenized tensor $\overline{A}:=\cH A$ coincides with the homogenization of $\widetilde{A}$, i.e. 
\[
\cH\left(\iota^{m\epsilon}\right)^{-1}\cD A^{\epsilon}=\cH\widetilde{A}=\cD\overline{A}=\cD\cH\left(\iota^{\epsilon}\right)^{-1}A^{\epsilon}.
\]
\end{lem}
\begin{proof}
To begin with, we apply the homogenization formula in Def. \ref{eq:homogenization-operator-def} to $\widetilde{A}$:
\begin{equation}
\cH\widetilde{A}\left(x\right)=\int_{\bT^{d}}\widetilde{A}\left(x,\lambda\right)+\widetilde{A}\left(x,\lambda\right)\nabla\widetilde{\chi}\left(\lambda;x\right)^{T}\d{\lambda}\label{eq:interchangability-exposition}
\end{equation}
where $\widetilde{\chi}$ solves the cell problem $-\nabla_{\lambda}\cdot\left(\nabla_{\lambda}\widetilde{\chi}\ \widetilde{A}^{T}\right)=\nabla_{\lambda}\cdot\widetilde{A}^{T}$. Since $\widetilde{A}\left(x,\lambda\right)=A\left(\phi\left(x\right),\lambda+\varphi_{k}\right)$ in $T_{k}$ and the cell problem does not contain partial differentials in $x$, one can easily verify that the solution reads $\widetilde{\chi}\left(\lambda;x\right)=\chi\left(\lambda+\varphi_{k};\phi\left(x\right)\right)$. Thus, with a change in the variable $\widetilde{\lambda}:=\lambda+\varphi_{k}$, Eqn. \ref{eq:interchangability-exposition} leads to%
\begin{comment}
\begin{align*}
\cH\widetilde{A}\left(x\right) & =\int_{\bT^{d}}A\left(\phi\left(x\right),\lambda+\varphi_{k}\right)+A\left(\phi\left(x\right),\lambda+\varphi_{k}\right)\nabla\chi\left(\lambda+\varphi_{k};\phi\left(x\right)\right)^{T}\d{\lambda}\\
 & =\int_{\bT^{d}}A\left(\phi\left(x\right),\widetilde{\lambda}\right)+A\left(\phi\left(x\right),\widetilde{\lambda}\right)\nabla\chi\left(\widetilde{\lambda};\phi\left(x\right)\right)^{T}\d{\widetilde{\lambda}}\\
 & =\overline{A}\left(\phi\left(x\right)\right)\\
 & =\cD\overline{A}\left(x\right).
\end{align*}
\end{comment}
\[
\cH\widetilde{A}\left(x\right)=\int_{\bT^{d}}A\left(\phi\left(x\right),\widetilde{\lambda}\right)+A\left(\phi\left(x\right),\widetilde{\lambda}\right)\nabla\chi\left(\widetilde{\lambda};\phi\left(x\right)\right)^{T}\d{\widetilde{\lambda}}=\cD\overline{A}\left(x\right).
\]
\end{proof}
\begin{rem}
\label{rem:commutative diagram}We conclude Lem. \ref{lem:interchangability} as the following commutative diagram.

\[\begin{tikzcd}[ampersand replacement=\&]
A^{\epsilon}
\arrow[d, "\cD"']
\arrow[r, "\left(\iota^{\epsilon}\right)^{-1}"]
\&
A
\arrow[r, "\cH"]
\&
\overline{A}
\arrow[d, "\cD"]
\\
\cD A^{\epsilon}
\arrow[r, "\left(\iota^{m\epsilon}\right)^{-1}"']
\&
\widetilde{A}
\arrow[r, "\cH"']
\&
\cD \overline{A}
\end{tikzcd}\]
\end{rem}

\subsubsection{Error decomposition}

So far, we have defined $u_{D}^{\epsilon}$ as the solution to the dilated oscillatory problem (Eqn. \ref{eq:dilation-solver-dilated-pde}) and $u_{0}$ as the solution to the homogenized problem (Eqn. \ref{eq:intro-homogenized-pde}). We may also introduce $u_{D,h}^{\epsilon}$ as the numerical solution to the dilated system (Eqn. \ref{eq:dilation-solver-dilated-pde}) where the mesh size is specified by $h$. In this section, we aim to build a set of error estimates on $u_{D,h}^{\epsilon}-u_{0}$ and find out how it depends on the parameters $L,m,$ and $\nu$. To this end, we can utilize the solution $u_{0,D}$ to Eqn. \ref{eq:dilation-solver-homogenized-dilated-pde}. By the triangle inequality,
\begin{equation}
\norm{u_{D,h}^{\epsilon}-u_{0}}\le\underbrace{\norm{u_{D,h}^{\epsilon}-u_{D}^{\epsilon}}}_{\text{discretization error}}+\underbrace{\norm{u_{D}^{\epsilon}-u_{0,D}}}_{\text{homogenization error}}+\underbrace{\norm{u_{0,D}-u_{0}}}_{\text{dilation error}}.\label{eq:full-error-estimate-tri}
\end{equation}
The homogenization error is treated similarly as in the section on partial dilation. According to Thm. \ref{thm:homogenization-error-Lp-estimate}, there exist constants $C\left(\Omega,r,M,d\right)$ and $r_{0}\left(\Omega,r,M,d\right)$ s.t. 
\begin{equation}
\norm{u_{D}^{\epsilon}-u_{0,D}}_{L^{2}\left(\Omega\right)}\le Cm\epsilon\ln\left(\frac{r_{0}}{m\epsilon}\right)\norm f_{L^{2}\left(\Omega\right)}.\label{eq:dilated-system-homogenization-error}
\end{equation}
The dilation error is addressed in Sec. \ref{subsec:Perturbation-analysis} by using a standard perturbation analysis approach. The discretization error of the FEM solution is analyzed in Sec. \ref{subsec:Discretization-error-fem}.

\subsubsection{Dilation error and regularity\label{subsec:Perturbation-analysis}}

According to the perturbation analysis \cite[Theorem 2.1]{Bonito2013}, the dilation error depends on the $L^{\infty}$-norm difference of the homogenized coefficients, i.e.

\[
\norm{u_{0}-u_{0,D}}_{H_{0}^{1}\left(\Omega\right)}\apprle\norm{\nabla u_{0}}_{L^{2}\left(\Omega\right)}\norm{\overline{A}-\cH\widetilde{A}}_{L^{\infty}\left(\Omega\right)}.
\]
Since $\left(\cH\widetilde{A}\right)\left(x\right)=\left(\cD\overline{A}\right)\left(x\right)=\overline{A}\left(\phi\left(x\right)\right)$ and $\phi\left(x\right)-x=\cO\left(L\right)$ (Lem. \ref{lem:shrinkage-difference-estimate}), the $L^{\infty}$-norm difference can be easily derived if $\overline{A}$ is smooth. To show the regularity of $\overline{A}$, we need to prepare a few key concepts and lemmas.
\begin{defn}
\label{def:spectral norm}For a matrix $B\in\bR^{d\times d}$, let $\norm B_{2}:=\sup_{\xi\neq0}\frac{\abs{B\xi}}{\abs{\xi}}$ denote the matrix norm induced by the vector 2-norm and $\norm B_{\max}:=\max_{i,j}\abs{B_{ij}}$ denote the max norm. The $L^{p}$ norm of a matrix-valued function $\cB\left(x\right)$ is defined as the $L^{p}$ norm of the point-wise spectral norm, i.e. $\norm{\cB}_{L^{p}\left(\Omega\right)}:=\norm{\norm{\cB\left(x\right)}_{2}}_{L^{p}\left(\Omega\right)}$.
\end{defn}

We now derive the regularity of the homogenized tensor $\overline{A}$.
\begin{lem}
\label{lem:HA regularity}We assume that a given $A\in\cS$ satisfies 
\begin{itemize}
\item $A\in C^{1}\left(\Omega\times\bT^{d};M^{d\times d}\right)$ and $A\left(x\right)$ is symmetric,
\item there exists $r>0$ s.t. $r\le\lambda_{\min}\left(A\left(x\right)\right)$ holds for any $x\in\Omega$, and
\item there exists a uniform bound $Q>0$ on $A_{ij}$, $\nabla_{x}A_{ij}$, $\nabla_{\lambda}A_{ij}$, $\nabla_{x}\nabla_{\lambda}A_{ij}$ for $1\le i,j\le d$ and $\left(x,\lambda\right)\in\Omega\times\bT^{d}$.
\end{itemize}
Then $\overline{A}=\cH A$ is Lipschitz in terms of the spectral norm, i.e. there exists a constant $C\left(\Omega,r,Q,d\right)$ s.t.
\[
\norm{\overline{A}\left(x\right)-\overline{A}\left(\widehat{x}\right)}_{2}\le C\left(\Omega,r,Q,d\right)\abs{x-\widehat{x}}
\]
holds for any $x,\widehat{x}\in\Omega$.
\end{lem}
\begin{proof}
We first need to establish a few estimates on the first order corrector field $\chi\left(\cdot;x\right)$. Let $a_{l}$ denote the $l$-th column of $A$. Since $\chi\left(\cdot;x\right),\chi\left(\cdot;\widehat{x}\right)\in\overline{H}^{1}\left(\bT^{d}\right)$ solve the cell problems
\begin{align*}
-\nabla_{\lambda}\cdot\left(A\left(x,\lambda\right)\nabla_{\lambda}\chi_{l}\left(\lambda;x\right)\right) & =\nabla_{\lambda}\cdot a_{l}\left(x,\lambda\right),\\
-\nabla_{\lambda}\cdot\left(A\left(\widehat{x},\lambda\right)\nabla_{\lambda}\chi_{l}\left(\lambda;\widehat{x}\right)\right) & =\nabla_{\lambda}\cdot a_{l}\left(\widehat{x},\lambda\right),
\end{align*}
by the perturbation theory for PDEs with periodic boundary conditions  (as an extension of \cite[Theorem 2.1]{Bonito2013}), the difference can be estimated by
\begin{multline}
\norm{\chi_{l}\left(\lambda;x\right)-\chi_{l}\left(\lambda;\widehat{x}\right)}_{\overline{H}^{1}\left(\bT^{d}\right)}\le\frac{1}{r}\norm{\nabla_{\lambda}\cdot\left(a_{l}\left(x,\lambda\right)-a_{l}\left(\widehat{x},\lambda\right)\right)}_{H^{-1}\left(\bT^{d}\right)}+\\
+\frac{1}{r}\norm{\nabla_{\lambda}\chi_{l}\left(\lambda;x\right)}_{L^{2}\left(\bT^{d}\right)}\norm{A\left(x,\lambda\right)-A\left(\widehat{x},\lambda\right)}_{L^{\infty}\left(\bT^{d}\right)}\label{eq:dilation-chi-property-perturbation}
\end{multline}
where the $\overline{H}^{1}\left(\Pi\right)$ semi-norm is defined as $\abs v_{\overline{H}^{1}\left(\Pi\right)}:=\norm{\nabla v}_{L^{2}\left(\Pi\right)}$. For the first term, we apply mean value theorem and bound the partial derivative
\begin{equation}
\norm{\nabla_{\lambda}\cdot\left(a_{l}\left(x,\lambda\right)-a_{l}\left(\widehat{x},\lambda\right)\right)}_{L^{2}\left(\bT^{d}\right)}\le\sum_{j=1}^{d}\left[\abs{\partial_{\lambda_{j}}\partial_{x}A_{jl}\left(\xi_{j},\lambda\right)}\abs{x_{j}-\widehat{x}_{j}}\right]\le dQ\abs{x-\widehat{x}}.\label{eq:dilation-chi-property-da}
\end{equation}
For the term involving $\chi$, we can apply a standard estimate in the Lax-Milgram theorem to conclude that
\begin{equation}
\norm{\nabla_{\lambda}\chi_{l}\left(\lambda;x\right)}_{L^{2}\left(\bT^{d}\right)}\le\frac{1}{r}\norm{\nabla_{\lambda}\cdot a_{l}\left(x,\lambda\right)}_{\overline{H}^{-1}\left(\bT^{d}\right)}\le\frac{C_{2}Q\sqrt{d}}{r}.\label{eq:dilation-chi-property-dchi}
\end{equation}
For the difference term in $A$ , we simply utilize the fact that each entry of $A$ is $C^{1}$ and that all matrix norms are equivalent:
\begin{equation}
\norm{A\left(x,\lambda\right)-A\left(\widehat{x},\lambda\right)}_{L^{\infty}\left(\bT^{d}\right)}\le C_{1}\left(d\right)dQ\abs{x-\widehat{x}}\label{eq:dilation-chi-property-diffA}
\end{equation}
where $C_{1}\left(d\right)$ is a constant s.t. $\frac{1}{C_{1}}\norm B_{2}\le\norm B_{\max}\le C_{1}\norm B_{2}$ holds for any matrix $B\in\bR^{d\times d}$.%
{} We combine Eqn. \ref{eq:dilation-chi-property-perturbation}, \ref{eq:dilation-chi-property-da}, \ref{eq:dilation-chi-property-dchi}, and \ref{eq:dilation-chi-property-diffA} to argue that
\begin{equation}
\norm{\chi_{l}\left(\lambda;x\right)-\chi_{l}\left(\lambda;\widehat{x}\right)}_{\overline{H}^{1}\left(\bT^{d}\right)}\le C_{\chi}\abs{x-\widehat{x}},C_{\chi}:=\frac{C_{2}dQ}{r}+\frac{C_{1}C_{2}d^{3/2}Q^{2}}{r^{2}},\label{eq:dilation-chi-property-Hbar1}
\end{equation}
which implies that $\chi_{l}$ is a Lipschitz function w.r.t. the parameter $x$ in terms of the $\overline{H}^{1}\left(\bT^{d}\right)$ norm.

Now we move on to the homogenized coefficient. %
Due to the fact that $\norm{\cdot}_{2}$ is sub-multiplicative and $\norm{\int_{\bT^{d}}B\d{\lambda}}_{2}\le\int_{\bT^{d}}\norm B_{2}\d{\lambda}$, it suffices to estimate the pointwise difference between $A\left(x,\cdot\right)$ and $A\left(\widehat{x},\cdot\right)$. Using the Lipschitz continuity of $A$, we have
\begin{equation}
\norm{A\left(x,\lambda\right)-A\left(\widehat{x},\lambda\right)}_{2}\le C_{1}\max_{i,j}\abs{\sup_{\xi}\nabla_{x}A_{ij}\left(\xi,\lambda\right)}\abs{x-\widehat{x}}\le C_{1}Q\abs{x-\widehat{x}}\label{eq:dilation-homogenized-I1}
\end{equation}
which implies that
\begin{equation}
\norm{A\left(x,\lambda\right)-A\left(\widehat{x},\lambda\right)}_{2}\norm{\nabla_{\lambda}\chi\left(\lambda;x\right)}_{2}\le\frac{C_{1}^{2}C_{2}\sqrt{d}Q^{2}}{r}\abs{x-\widehat{x}}.\label{eq:dilation-homogenized-I2}
\end{equation}
Meanwhile, Eqn. \ref{eq:dilation-chi-property-Hbar1} implies 
\begin{equation}
\norm{A\left(\widehat{x},\lambda\right)}_{2}\norm{\nabla_{\lambda}\chi\left(\lambda;x\right)-\nabla_{\lambda}\chi\left(\lambda;\widehat{x}\right)}_{2}\le C_{1}QC_{\chi}\abs{x-\widehat{x}}.\label{eq:dilation-homogenized-I3}
\end{equation}
Thus, a combination of triangular inequality and Eqn. \ref{eq:homogenization-operator-def}, \ref{eq:dilation-homogenized-I1}, \ref{eq:dilation-homogenized-I2}, and \ref{eq:dilation-homogenized-I3} leads to
\[
\norm{\overline{A}\left(x\right)-\overline{A}\left(\widehat{x}\right)}_{2}\le C_{1}Q\left(1+C_{\chi}+C_{1}C_{2}\sqrt{d}Q/r\right)\abs{x-\widehat{x}}
\]
for any $x,\widehat{x}\in\Omega$.
\end{proof}
\begin{prop}
\label{prop:dilation error estimate}Given $A\in\cS$ that satisfies the conditions in Lem. \ref{lem:HA regularity}, there exists a constant $C\left(\Omega,r,Q,d\right)$ s.t. the dilation error is controlled by
\begin{equation}
\norm{u_{0,D}-u_{0}}_{L^{2}\left(\Omega\right)}\le C\left(\Omega,r,Q,d\right)L\norm{\nabla u_{0}}_{L^{2}\left(\Omega\right)}.\label{eq:dilation-error-estimate}
\end{equation}
\end{prop}
\begin{proof}
By the perturbation analysis (see \cite[Theorem 2.1]{Bonito2013}), we have
\begin{equation}
\norm{\nabla\left(u_{0,D}-u_{0}\right)}_{L^{2}\left(\Omega\right)}\le\frac{1}{r}\norm{\nabla u_{0}}_{L^{2}\left(\Omega\right)}\norm{\overline{A}-\cD\overline{A}}_{L^{\infty}\left(\Omega\right)}.\label{eq:nabla error bound}
\end{equation}
By definition of the dilation operator and the estimate in Lem. \ref{lem:HA regularity}, we have%
\begin{comment}
\begin{align*}
\norm{\overline{A}-\cD\overline{A}}_{L^{\infty}\left(\Omega\right)} & =\max_{x\in\Omega}\norm{\overline{A}\left(x\right)-\cD\overline{A}\left(x\right)}_{2}\\
 & =\max_{x\in\Omega}\norm{\overline{A}\left(x\right)-\overline{A}\left(\phi\left(x\right)\right)}_{2}\\
 & \le\max_{x\in\Omega}C\left(\Omega,r,Q,d\right)\abs{x-\phi\left(x\right)}\\
 & \le C\left(\Omega,r,Q,d\right)dL.
\end{align*}
\end{comment}
\begin{equation}
\norm{\overline{A}-\cD\overline{A}}_{L^{\infty}\left(\Omega\right)}=\max_{x\in\Omega}\norm{\overline{A}\left(x\right)-\overline{A}\left(\phi\left(x\right)\right)}_{2}\le C\left(\Omega,r,Q,d\right)dL.\label{eq:ADA estimate}
\end{equation}
Then, Eqn. \ref{eq:dilation-error-estimate} is established by combining Eqn. \ref{eq:nabla error bound} and \ref{eq:ADA estimate} and applying the Poincare inequality to $u_{0,D}-u_{0}$.
\end{proof}
\begin{rem*}
We have shown that the dilation error scales linearly with the mesoscopic length $L$ in Prop. \ref{prop:dilation error estimate}. Nevertheless, this bound is not optimal since a faster convergence is observed for $\nu=1/2$ from the numerical experiments; we refer the readers to Sec. \ref{par:numerical-dilation-error} for more details.
\end{rem*}

\subsection{Numerical aspects\label{subsec:numerical-aspects-and-concerns}}

\subsubsection{Finite element approach\label{subsec:Finite-element-approach}}

The finite element method is applied to these problems. We use the PDE Toolbox in MATLAB to generate meshes and solve for the solution. It is worth pointing out that there is a natural restriction on the parameters. First, the mesh size shall resolve the effective parameter $\widetilde{\epsilon}=m\epsilon$, i.e. $h\ll m\epsilon$. Meanwhile, the mesoscopic length $L$ usually needs to accommodate at least one oscillation before the dilation operator ``resets'', which implies $m\epsilon<L$. We call $h\ll m\epsilon<L$ as the scale condition.

\subsubsection{A priori estimate on discretization error (1D case)\label{subsec:Discretization-error-fem}}

In general, the dilated coefficient $\cD A^{\epsilon}$ has discontinuities on the mesoscopic scale $L$ and thus the classical FEM/FDM error estimate no longer applies. Non-continuous coefficients are usually a challenge, but for the case of one dimensional problems, we can utilize the closed form of the solution in an integral form. In fact, we can show a first order convergence with respect to the mesh size $h$. We study the second order problem 
\begin{equation}
-\left(a\,u'\right)'=f,x\in\Omega=\left[0,1\right]\label{eq:1d-proof-strong}
\end{equation}
under the following assumption.
\begin{assumption}
The domain $\Omega$ is partitioned by an evenly-spaced mesh $\left\{ x_{i}\right\} $, $x_{i}:=i\,h$ with mesh size $h$ and an evenly-spaced meso-mesh $\left\{ X_{j}\right\} ,X_{j}:=j\,L$ with length $L$; we assume $1/h,1/L\in\bZ^{+}$. The coefficient function $a$ is piecewise-$C^{1}$ on $\Omega_{j}:=\left[X_{j-1},X_{j}\right]$ for $j=1,\dots,1/L$. Moreover, there exist $M,\epsilon>0$ s.t. $M^{-1}\le a\le M,\ \abs{a'}\le M\epsilon^{-1}$ for every $x\in\Omega_{j}$ and all $j$. We call $\left(h,L\right)$ an aligned-configuration if $L/h\in\bZ^{+}$ or a non-aligned configuration otherwise.
\end{assumption}
We introduce the variational problem and the numerical counterpart as follows. Then, the numerical error is related to the interpolation error through the Cea's lemma.
\begin{defn}
Eqn. \ref{eq:1d-proof-strong} is understood in the following variational sense: find $u\in H_{0}^{1}\left(\Omega\right)$s.t. $\forall v\in H_{0}^{1}\left(\Omega\right)$,
\begin{equation}
\left\langle a\,u',v'\right\rangle _{L^{2}\left(\Omega\right)}=\left(f,v\right)_{H^{-1}\left(\Omega\right),H^{1}\left(\Omega\right)}\label{eq:1d-proof-variational}
\end{equation}
where $\left\langle \cdot,\cdot\right\rangle _{L^{2}}$ is the $L^{2}$-inner product and $\left(\cdot,\cdot\right)_{H^{-1},H^{1}}$ is the duality pairing. Given the piece-wise linear interpolation operator $I_{h}$ and $V_{h}:=I_{h}H_{0}^{1}\left(\Omega\right)$, the numerical solution $u_{h}\in V_{h}$ solves Eqn. \ref{eq:1d-proof-variational} with any $v_{h}\in V_{h}$.
\end{defn}

Since $u'$ is piecewise $H^{1}$ within each $\Omega_{j}$, we can obtain an $H^{2}$-norm estimate in an aggregated manner by integrating $u''$ on each sub-interval.
\begin{prop}
\label{prop:H2-aggregated-norm}Given $u$ solving Eqn. \ref{eq:1d-proof-variational}, the aggregated $H^{2}$-norm is bounded by
\[
\left|u\right|_{H^{2}}^{\dagger}:=\sqrt{\sum_{j=1}^{L^{-1}}\norm{u''}_{L^{2}\left(\Omega_{j}\right)}^{2}}\le\left(M+M^{3}\epsilon^{-1}+M^{5}\epsilon^{-1}\right)\norm f_{L^{2}}.
\]
Besides, the $L^{\infty}$ bounds are given by $\norm u_{L^{\infty}},\norm{u'}_{L^{\infty}}\le M\norm f_{L^{1}}+M^{3}\norm f_{L^{2}}$.
\end{prop}
\begin{proof}
We write the solution $u$ in its closed form $u\left(x\right)=\int_{0}^{x}a^{-1}\left(y\right)\left[F\left(y\right)+C\right]\d y$ where $F\left(x\right):=\int_{0}^{x}f\left(y\right)\d y$ and $C:=-\int_{0}^{1}\frac{F}{a}\d y/\int_{0}^{1}\frac{1}{a}\d y$. The second order derivative can be written explicitly as $u''=\frac{f}{a}-\left(F+C\right)\frac{a'}{a^{2}}$ for $x\in\Omega_{j}$. Since $F\left(x\right)$ vanishes at $x=0$, by Poincare inequality, $\norm F_{L^{2}}\le\norm f_{L^{2}}$ and subsequently $\norm{F\frac{a'}{a^{2}}}_{L^{2}}\le M^{3}\epsilon^{-1}\norm f_{L^{2}}$. Similarly, $\abs C\le M^{2}\norm F_{L^{1}}\le M^{2}\norm f_{L^{2}}$. Thus, the $H^{2}$-semi-norm can be bounded by
\[
\sum_{j=1}^{L^{-1}}\norm{u''}_{L^{2}\left(\Omega_{j}\right)}^{2}\le\left(M+M^{3}\epsilon^{-1}+M^{5}\epsilon^{-1}\right)^{2}\norm f_{L^{2}}^{2}.
\]
The $L^{\infty}$ bounds are derived by $\abs{u'}\le\abs{Fa^{-1}}+\abs C\abs{a^{-1}}\le M\norm f_{L^{1}}+M^{3}\norm f_{L^{2}}$ and $\abs{u\left(x\right)}\le\max_{\xi\in\left[0,x\right]}\abs{u'\left(\xi\right)}$.
\end{proof}
\begin{prop}
\label{prop:approximation-error}The approximation error is $\norm{u-I_{h}u}_{H^{1}}\le\sqrt{h^{4}+h^{2}}\left|u\right|_{H^{2}}^{\dagger}$ for aligned configurations. 
\end{prop}
\begin{proof}
The interpolation difference $\delta u:=u-I_{h}u$ vanishes on the boundary of each interval $\Delta_{i}:=\left(x_{i-1},x_{i}\right)$, thus by Poincare inequality
\begin{equation}
\norm{\left(u-I_{h}u\right)'}_{L^{2}\left(\Delta_{i}\right)}\le h\left|\left(u-I_{h}u\right)'\right|_{H^{1}\left(\Delta_{i}\right)}=h\left|u\right|_{H^{2}\left(\Delta_{i}\right)}\label{eq:H1 interpolation error}
\end{equation}
and subsequently
\begin{equation}
\norm{u-I_{h}u}_{L^{2}\left(\Delta_{i}\right)}\le h\left|u-I_{h}u\right|_{H^{1}\left(\Delta_{i}\right)}=h^{2}\left|u\right|_{H^{2}\left(\Delta_{i}\right)}.\label{eq:L2 interpolation error}
\end{equation}
Since $\Delta_{i}$ is always included in some $\Omega_{j}$ for aligned configurations, we have
\[
\sum_{i=1}^{h^{-1}}\norm{u-I_{h}u}_{H^{1}\left(\Delta_{i}\right)}^{2}\le\left(h^{4}+h^{2}\right)\sum_{i=1}^{h^{-1}}\left|u\right|_{H^{2}\left(\Delta_{i}\right)}^{2}\le\left(h^{4}+h^{2}\right)\left(\left|u\right|_{H^{2}}^{\dagger}\right)^{2}.
\]
\end{proof}
\sepia{Cea's lemma \cite[Theorem 12.3]{quarteroniNumericalMathematics2007} states that $\norm{u-u_{h}}_{H^{1}\left(\Omega\right)}\le2M^{2}\norm{u-I_{h}u}_{H^{1}\left(\Omega\right)}$ which helps to bound the error in solution by the interpolation error. Consequently, we obtain the $H^{1}$ and $L^{2}$ estimates by combining this with Prop. \ref{prop:H2-aggregated-norm} and \ref{prop:approximation-error}.}

\begin{thm}
\label{thm:L2-estimate}Given $h<\epsilon<L<1<M$, the $H^{1}$-error is $\norm{u-u_{h}}_{H^{1}}\le6\sqrt{2}M^{7}\epsilon^{-1}h\norm f_{L^{2}}$ and the $L^{2}$-error reads $\norm{u-u_{h}}_{L^{2}}\le36M^{13}\epsilon^{-2}h^{2}\norm f_{L^{2}}$. 
\end{thm}
\begin{proof}
We apply the Aubin-Nitsche duality method. Let $e:=u-u_{h}$ and we seek for solution $\phi$ to $\left\langle a\,\phi',\psi'\right\rangle =\left(e,\psi\right)$, $\forall\psi\in H_{0}^{1}\left(\Omega\right)$. Since $e$ is the residual, it is orthogonal to $w_{h}\in V_{h}$ in terms of the $a$-inner product, leading to
\[
\norm e_{L^{2}}^{2}=\left\langle a\,\phi',e'\right\rangle =\left\langle a\,e',\left(\phi-I_{h}\phi\right)'\right\rangle \le M\norm e_{H^{1}}\norm{\phi-I_{h}\phi}_{H^{1}}.
\]
By \ref{prop:approximation-error}, the interpolation error is bounded by $\norm{\phi-I_{h}\phi}_{H^{1}}\le\sqrt{2}h\left|\phi\right|_{H^{2}}^{\dagger}$ where the regularity is established in Prop. \ref{prop:H2-aggregated-norm} as $\left|\phi\right|_{H^{2}}^{\dagger}\le3M^{5}\epsilon^{-1}\norm e_{L^{2}}$ and $\norm{\phi}_{L^{\infty}},\norm{\phi'}_{L^{\infty}}\le2M^{3}\norm e_{L^{2}}$. Thus, we have $\norm e_{L^{2}}\le3\sqrt{2}hM^{6}\epsilon^{-1}\norm e_{H^{1}}$.
\end{proof}
\begin{rem*}
For non-aligned configurations, the asymptotic $L^{2}$-bound is of order $\cO\left(h\right)$, which is slower than the aligned counterpart. In practice, the alignment is easy to satisfy since one can pick any $L$ and set $h$ as any fraction of $L$ as long as they meet the scale condition (see Sec. \ref{subsec:Finite-element-approach}).
\end{rem*}
\begin{rem*}
In general, for dimension $d>1$, it is much more difficult to build analogies to the approximation result \ref{prop:approximation-error} since $u$ may not be in $H^{2}\left(\Omega_{j}\right)$ on the sub-domains. Another technical issue is posed by irregularity near vertices that are shared by neighboring elements. In fact, it is pointed out in \cite{ciarletDomainDecompositionMethods2017} that there exists an exponent $r_{\max}\in\left(0,1\right]$ s.t. $u\in H^{1+r}\left(\Omega_{j}\right)$ for $0\le r<r_{\max}$ and $\norm u_{H^{1+r}}\apprle\norm f_{L^{2}}$. A domain-decomposition based method is proposed to handle the missing regularity for the mixed problem.
\end{rem*}
We can now put the pieces together (Eqn. \ref{eq:full-error-estimate-tri} and \ref{eq:dilated-system-homogenization-error}, Prop. \ref{prop:dilation error estimate}, Cor. \ref{thm:L2-estimate}) and derive an error bound for the FEM solution to the dilated system.
\begin{thm}
\label{cor:fem-solver-full-error}For aligned 1d domain setups and assumptions made in Assump. \ref{assu:conditions on A}, Cor. \ref{thm:L2-estimate}, and Prop. \ref{prop:dilation error estimate}, there exists constants $r_{0}\left(\Omega,r,M,d\right)$ and $C$ depending on $\Omega,r,Q,d,M$ s.t. the full error is bounded by
\[
\norm{u_{D,h}^{\epsilon}-u_{0}}_{L^{2}\left(\Omega\right)}\le C\left[\epsilon^{-2}h^{2}+m\epsilon\ln\left(\frac{r_{0}}{m\epsilon}\right)+L\right].
\]
\end{thm}

\subsection{Numerical experiments\label{subsec:experiments}}

We present the following two numerical examples and verify the asymptotic orders of each error component. We also include an integrated test to compare the performance of methods of partial dilation, with or without scale identification, and local dilation.

\subsubsection{Examples set-up}

According to Eqn. \ref{eq:homogenization-operator-def} and \ref{eq:intro-homogenized-pde}, the homogenized solution $u$ can be solved as long as the homogenized tensor $\overline{A}$ is obtained. However, it involves solving $\chi$ from Eqn. \ref{eq:intro-cell-problem} which in general is a family of elliptic PDEs, rendering this approach computationally intractable in higher dimensions $d>1$. Nevertheless, a closed form solution can be obtained if the coefficient only depends on one of the fast coordinates, which is often referred to as the case of layered materials. Without loss of generality, we work with underlying dimension $d=2$. 

\paragraph{Layered material \label{par:Layered-material}}

We start with a case where a closed form for $\overline{A}$ is available (\cite[Sec. 12.6.2]{Grigorios2008}). For a given diffusion coefficient in the form of $A\left(x,\lambda\right)=A\left(x_{1},x_{2},\lambda_{1}\right)$, the homogenized tensor reads $\overline{A}_{11}=1/\left\langle 1/A_{11}\right\rangle $, $\overline{A}_{12}=\left\langle A_{12}/A_{11}\right\rangle /\left\langle 1/A_{11}\right\rangle $, and $\overline{A}_{22}=\overline{A_{12}}\overline{A_{21}}/\overline{A_{11}}-\left\langle A_{22}-A_{12}A_{21}/A_{11}\right\rangle $ where $\left\langle g\right\rangle \left(x_{1},x_{2}\right):=\int_{0}^{1}g\left(x_{1},x_{2},\lambda_{1}\right)\d{\lambda_{1}}$ is the averaging operation. 

For given tilting angle $\theta$ and intensity factor $\eta\in\left[0,1\right]$, let $A_{\text{lay}}^{\epsilon}$ denote 
\begin{equation}
A_{\text{lay}}^{\epsilon}\left(x_{1},x_{2};\eta,\theta\right):=\left[y_{1}\left(x_{1},x_{2}\right)+y_{2}\sin\left(\frac{2\pi}{\epsilon}\left(x_{1}\cos\theta-x_{2}\sin\theta\right)\right)\right]I_{2\times2}\label{eq:setting-layered 2cylinder smooth}
\end{equation}
where the diagonal entry reads $y_{1}\left(x_{1},x_{2}\right):=1+0.1x_{1}+0.05x_{2}$ and $y_{2}:=0.9\eta$. The source term reads
\[
f_{\text{lay}}\left(x_{1},x_{2}\right)=e^{-60\left(\left(x_{1}-0.3\right)^{2}+\left(x_{2}-0.3\right)^{2}\right)}+e^{-60\left(\left(x_{1}-0.7\right)^{2}+\left(x_{2}-0.7\right)^{2}\right)}.
\]
We shall point out that Eqn. \ref{eq:homogenization-operator-def} does not apply to $A_{\text{lay}}^{\epsilon}$ directly since the oscillation depends on the linear combination $X_{1}:=x_{1}\cos\theta-x_{2}\sin\theta$. A rotation is needed to obtain the correct homogenized tensor 
\[
\overline{A}_{\text{lay}}=U_{\theta}\overline{U_{\theta}^{T}A_{\text{lay}}^{\epsilon}U_{\theta}}U_{\theta}^{T},U_{\theta}:=\left(\begin{array}{cc}
\cos\theta & \sin\theta\\
-\sin\theta & \cos\theta
\end{array}\right).
\]

\paragraph{Material with heterogeneous oscillation \label{par:Material-with-heterogeneous}}

\begin{figure}[tbph]
\begin{centering}
\includegraphics[width=1\linewidth]{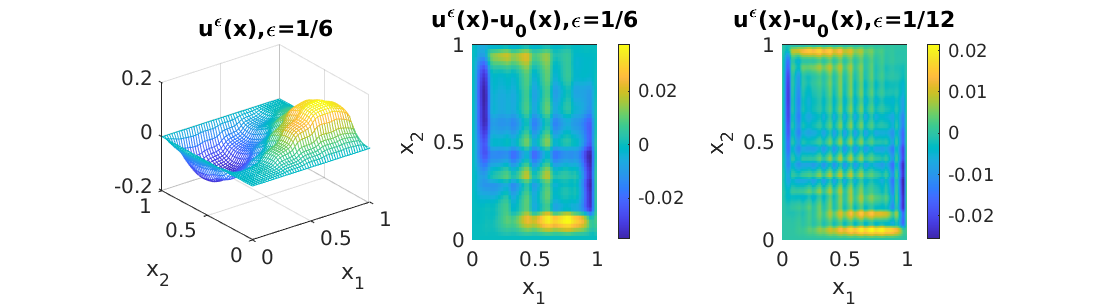}
\par\end{centering}
\caption{An illustrative demonstration of the oscillatory solution $u^{\epsilon}$ to the system described in Par. \ref{par:Material-with-heterogeneous} and Eqn. \ref{eq:setting-het osci2}. Left: the solution landscape when $\epsilon=\frac{1}{6}$. Middle/right: the residual $u^{\epsilon}-u_{0}$ for $\epsilon=\frac{1}{6}$ and $\frac{1}{12}$, respectively. Notice that the range on the right has shrunk by a half approximately compared to the one in the middle. The spatial frequency also doubled as $\epsilon$ is cut down in half.}

\label{fig:oscillatory solution}
\end{figure}

The second case is composed of oscillations that involves $x_{1}/\epsilon$ and $x_{2}/\epsilon$ and changing amplitudes from site to site. We define 
\begin{equation}
A_{\text{het}}^{\epsilon}\left(x_{1},x_{2};\eta\right)=\begin{bmatrix}1+\frac{1}{10}x_{1}+\frac{4}{5}\eta\left(\frac{4+x_{1}}{5}\right)s_{1} & \frac{1}{10}\eta\left(1-\frac{3}{10}x_{2}\right)s_{1}\\
\frac{1}{10}\eta\left(1-\frac{3}{10}x_{2}\right)s_{1} & 1+\frac{4}{5}\eta\left(\frac{7+3\sin\left(2\pi x_{2}\right)}{10}\right)s_{2}
\end{bmatrix}.\label{eq:setting-het osci2}
\end{equation}
The source term reads $f_{\text{het}}\left(x_{1},x_{2}\right)=10\left(x_{1}-x_{2}\right)$. The challenge lies in that for each evaluation of $\overline{A}_{\text{het}}\left(\widehat{x}_{1},\widehat{x}_{2}\right)$, one has to solve the parametrized cell problem and the cost can soon explode since the number of queries $\left\{ \left(\widehat{x}_{1},\widehat{x}_{2}\right)\right\} $ scales with the mesh size $h$ by $\cO\left(h^{-d}\right)$. One way to circumvent this situation is to exploit the regularity of $\overline{A}_{\text{het}}$ according to Lem. \ref{lem:HA regularity}. The cell problems are solved on a much coarser mesh and the homogenized tensor is interpolated at the query points. A figure containing solutions for finite $\epsilon$ is attached in Fig. \ref{fig:oscillatory solution}. The oscillatory solution $u^{\epsilon}$ can be decomposed into a smooth part $u_{0}$ and the first-order oscillation where the spatial frequency increases as $\epsilon\to0$.

\subsubsection{Empirical order of error components}

\paragraph{Dilation error $\protect\norm{u_{0,D}-u_{0}}_{L^{2}}$ \label{par:numerical-dilation-error}}

\begin{figure}[tbph]
\begin{centering}
\par\end{centering}
\begin{centering}
\par\end{centering}
\begin{centering}
\includegraphics[width=1\linewidth]{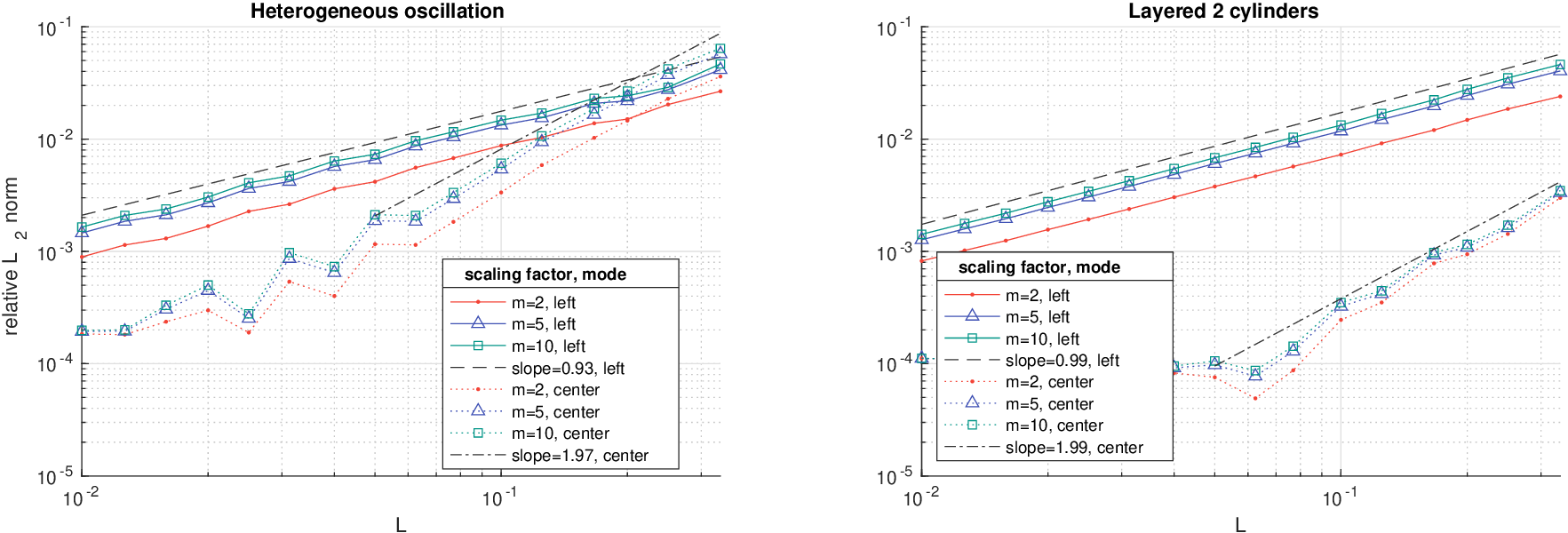}
\par\end{centering}
\caption{The dilation error scales with the mesoscopic scale $L$. A cross-comparison between different systems (described in Eqn. \ref{eq:setting-layered 2cylinder smooth} and \ref{eq:setting-het osci2}) is displayed from left to right. The mesh size is set to $h=0.15L$ adaptively.}
\label{fig:dilation-error}
\end{figure}

We start with the dilation error, shown in Fig. \ref{fig:dilation-error}. Predicted by Eqn. \ref{eq:dilation-error-estimate}, the relative error $\norm{u_{0,D}-u_{0}}_{L^{2}}/\norm{u_{0}}_{L^{2}}$ exhibits a linear convergence w.r.t. the mesoscopic scale $L$. We test and compare two type of dilations: anchoring to left ($\nu=0$) and to center ($\nu=0.5$). In general, the center dilation achieves a much better approximation and converges at a quadratic speed. Also, it is interesting to find out that a larger $m$ leads to a slightly larger error, which can be explained by the geometric nature of the dilation operator. With that being said, a larger $m$ implies a larger effective $\widetilde{\epsilon}=m\epsilon$ and thus cuts down the computational cost, so it is often a preferable trade-off.

\paragraph{Homogenization error $\protect\norm{u_{D}^{\epsilon}-u_{0,D}}_{L^{2}}$}

\begin{figure}[tbph]
\begin{centering}
\par\end{centering}
\begin{centering}
\includegraphics[width=1\linewidth]{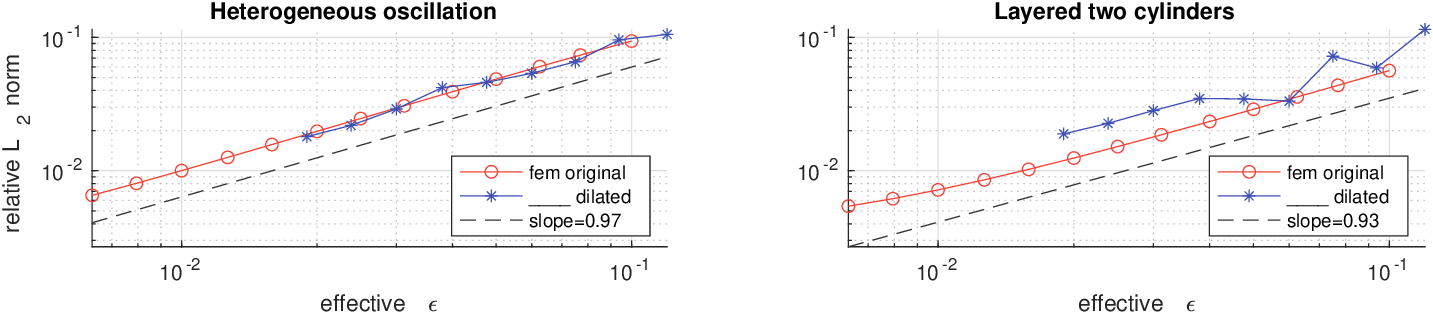}
\par\end{centering}
\caption{The homogenization error scales with the effective parameter $\widetilde{\epsilon}=m\epsilon$. Solutions are obtained by FEM with mesh size set to $h=0.1\widetilde{\epsilon}=0.1m\epsilon$. A cross-comparison between different systems (described in Eqn. \ref{eq:setting-layered 2cylinder smooth} and \ref{eq:setting-het osci2}) is displayed from left to right. The circles correspond to no dilation ($m=1$) and the star marks are obtained by dilating the coefficients with $L=0.1,m=3$. }

\label{fig:homogenization-error}
\end{figure}

We move onto the homogenization error, especially with the dilation in effect. We show the numerical results in Fig. \ref{fig:homogenization-error} with fitted asymptotic orders. Predicted by Thm. \ref{thm:homogenization-error-Lp-estimate}, we shall expect a first order convergence in the effective $\widetilde{\epsilon}$, considering the dilation effect. This convergence behavior is well verified by the numerical outcomes.

\paragraph{Discretization error $\protect\norm{u_{D,h}^{\epsilon}-u_{D}^{\epsilon}}_{L^{2}}$ \label{par:Discretization-error}}

\begin{figure}[tbph]
\begin{centering}
\par\end{centering}
\begin{centering}
\includegraphics[width=1\linewidth]{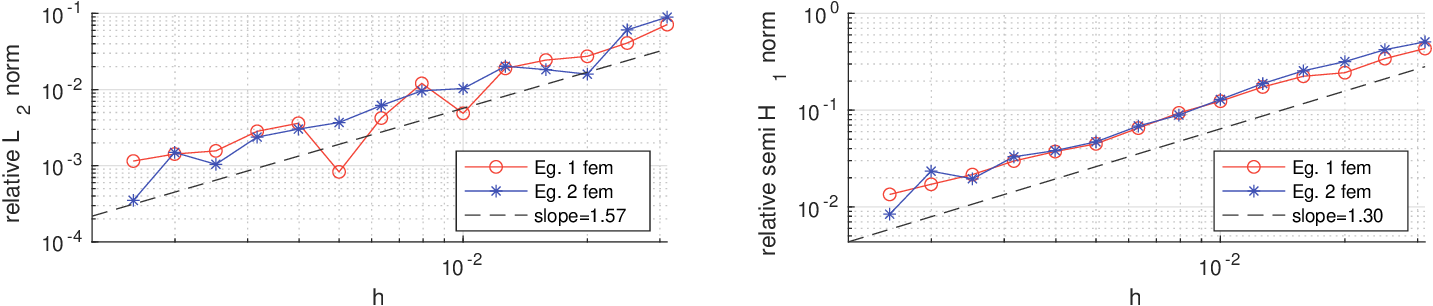}
\par\end{centering}
\caption{The discretization error scales with mesh size $h$ where the parameters are set to $\epsilon=0.04,m=2,L=0.1$. A cross-comparison between different systems (described in Eqn. \ref{eq:setting-layered 2cylinder smooth} and \ref{eq:setting-het osci2}) is displayed in different colors. Relative $L^{2}$ and semi-$H^{1}$ norms of the residual $u_{D,h}^{\epsilon}-u_{D}^{\epsilon}$ are examined in the left and right column respectively.}

\label{fig:discretization-error}
\end{figure}

We numerically verify the asymptotic order in Cor. \ref{thm:L2-estimate} and the results are shown in Fig. \ref{fig:discretization-error}. The $L^{2}$ and semi-$H^{1}$ norm of the residual $u_{D,h}^{\epsilon}-u_{D}^{\epsilon}$ decreases with mesh size $h$ and the order of convergence matches the expectation.

\subsubsection{Integrated test\label{subsec:integrated-test}}

\sepia{To conclude this section}, we demonstrate how well the dilation solver performs when it is applied as a drop-in replacement. Recall that our goal is to provide a middle ground where approximated solutions can be obtained by a relatively low cost. A direct way to cut computational complexity is to make the mesh coarser, i.e. to increase mesh size $h$. Judging by the scale condition $h\ll m\epsilon<L$, a larger $h$ requires a larger scaling factor $m$ and a larger mesoscopic length $L$. Since the scaling factor $m$ directly associates with $h$ and $L$ and further determines the computational cost, we pick it as the testing variable in the following experiment.

The test is set up as follows. A multiscale problem is given (as described in Eqn. \ref{eq:setting-layered 2cylinder smooth} and \ref{eq:setting-het osci2}) and the microscopic scale parameter $\epsilon$ is fixed. The scaling factor $m$ varies within a fixed testing range. To make the dilation solver a drop-in replacement, we only have access to $A^{\epsilon}$ rather than $A=\left(\iota^{\epsilon}\right)^{-1}A^{\epsilon}$. The other parameters, $h$ and $L$, are determined based on $m$ and a rough estimate on $\epsilon$. The dilated multiscale problem (Eqn. \ref{eq:dilation-solver-dilated-pde}) $-\nabla\cdot\left[\cD_{m,L}A^{\epsilon}\nabla u_{D}^{\epsilon}\right]=f$ and the homogenized problem $-\nabla\cdot\left[\cH A^{\epsilon}\nabla u_{0}\right]=f$ are solved and solutions are compared in the $L^{2}\left(\Omega\right)$ space. Furthermore, we also test the method of partial dilation where only the fast scale is dilated, i.e. $-\nabla\cdot\left[A^{m\epsilon}\nabla\widehat{u}^{\epsilon}\right]=f$. The error in the solution $\widehat{u}^{\epsilon}$ is expected to be $\cO\left(m\epsilon+h^{2}/\epsilon\right)$ as analyzed in Sec. \ref{subsec:Averaging-and-partial-dilation}. In the absence of the scale separated form, we solve $\left(\iota^{\epsilon}\right)^{-1}$ numerically by the following procedure. First, the observed signal $\left\{ A\left(x_{j}\right)\right\} $ is separated into the components with low and high frequencies by a mollifier kernel with proper support. The smooth (i.e. of low frequency) component is then projected on a monomial space for sparse identification. We deploy the EMD algorithm \cite{Huang1998} to extract the intrinsic mode functions from the high frequency residual. The dilation is thus only applied on each intrinsic modes.

\begin{figure}[tbph]
\begin{centering}
\par\end{centering}
\begin{centering}
\includegraphics[width=1\linewidth]{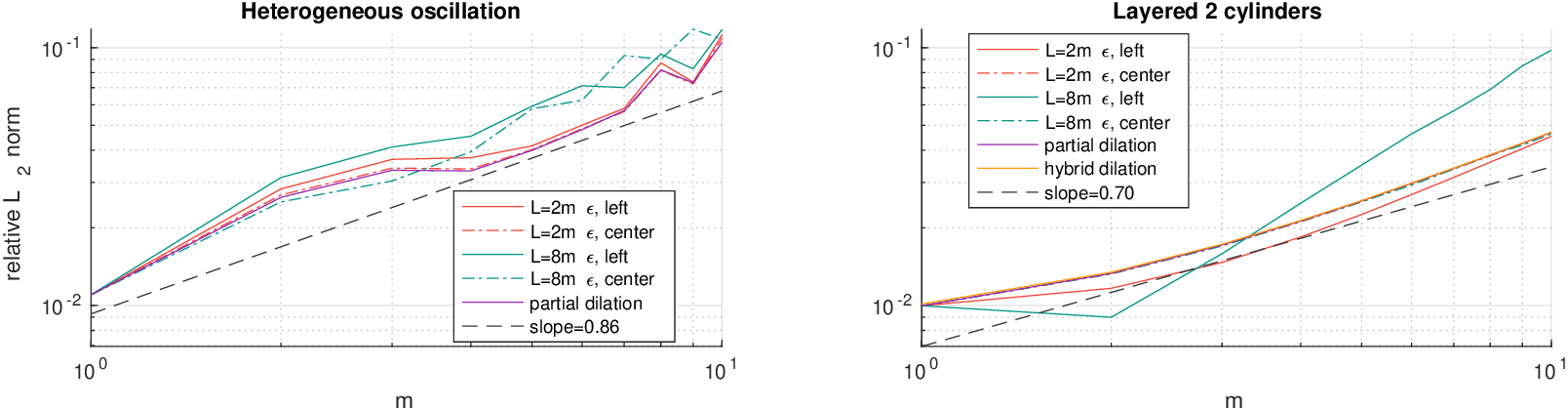}
\par\end{centering}
\caption{Error $\protect\norm{u_{D,h}^{\epsilon}-u_{0}}$ and $\protect\norm{\widehat{u}^{\epsilon}-u_{0}}$ plotted against the scaling factor $m$. A cross-comparison between different systems (described in Eqn. \ref{eq:setting-layered 2cylinder smooth} and \ref{eq:setting-het osci2}) is displayed from left to right. The $\epsilon$ parameter is set to $0.008$ and $m$ ranges from 1 to 10. The mesh size is $h=\frac{m\epsilon}{6.5}$ and the mesoscopic length is set to $L=2m\epsilon$ or $8m\epsilon$. We show the partial dilation approach (Eqn. \ref{eq:partial dilation-dilated-pde}) and the hybrid variant (marked as \textquotedblleft hybrid dilation\textquotedblright , described in Sec. \ref{subsec:Numerical-inversion}) especially in the layered material case. The smooth components are recovered exactly up to 4 digits of accuracy.}

\label{fig:integrated-tests}
\end{figure}

The error of different methods are plotted in Fig. \ref{fig:integrated-tests}. The error curves are quite close and they exhibit an almost first order relationship w.r.t. the scaling factor $m$. Since in these cases the mesh size $h$ is in proportion of $m$, a large $m$ certainly leads to a simpler mesh and a smaller linear system and thus it cuts down the computational complexity. The dilation method, regardless of the choice of $L$, performs as well as the partial dilation, showing that it is capable of relaxing the stiffness inside the problem without using the scale-separated form. Besides, the hybrid approach has an almost identical performance as the partial dilation for systems with a layered nature.

\subsection{Structure-aware dilation\label{subsec:Structure-aware-dilation}}

\begin{figure}[tbph]
\begin{centering}
\subfloat[Illustration on the channel effect. From left to right (showing the first component $\cdot_{11}$ only): original oscillatory tensor $A^{\epsilon}=A_{s}+A_{o}$, naively dilated tensor $\protect\cD A^{\epsilon}$, structure-aware dilated tensor $\protect\cD^{\text{s}}A=A_{s}+\protect\cD A_{o}$, and homogenized tensor $\overline{A}$. Notice that the channel is broken after applying naive dilation.]{\begin{centering}
\includegraphics[width=0.9\linewidth]{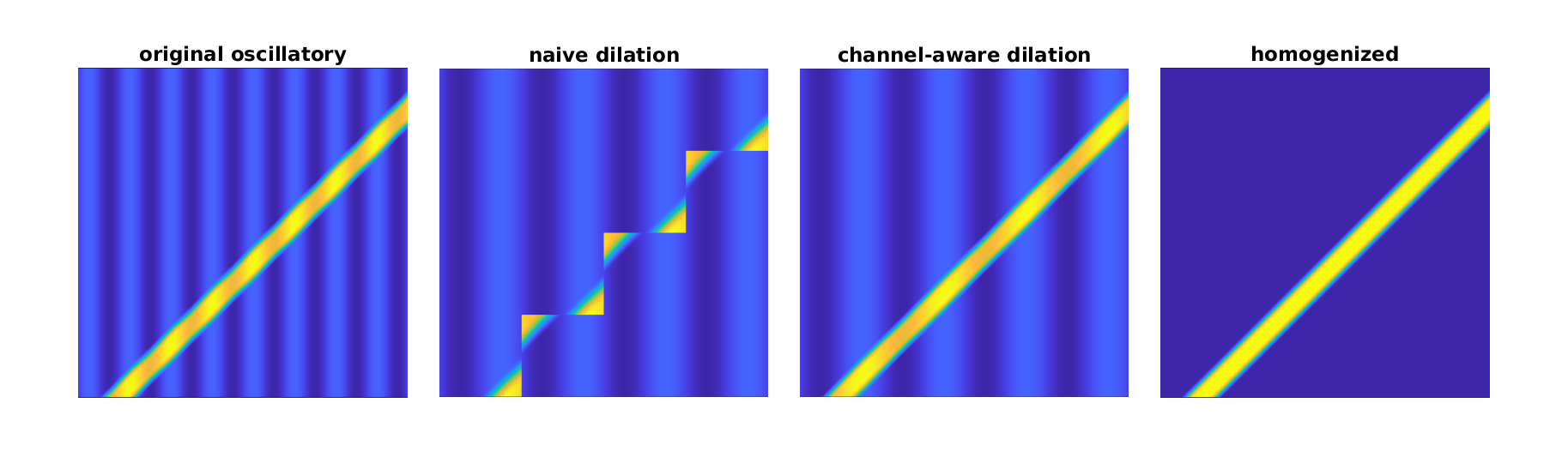}
\par\end{centering}
\label{fig:channel aware-dilation coeff}}
\par\end{centering}
\begin{centering}
\subfloat[Comparing the solution $u^{\epsilon}$ from the naive/structure-aware dilated system with the homogenized solution $u_{0}$. From left to right: solution to naively dilated tensor $\protect\cD A^{\epsilon}$, solution to structure-aware dilated tensor $\protect\cD^{\text{s}}A$, and homogenized solution to tensor $\overline{A}$. The \textquotedblleft +\textquotedblright{} sign indicates the center of the source while the cross sign is the center of the sink. ]{\begin{centering}
\includegraphics[width=0.9\linewidth]{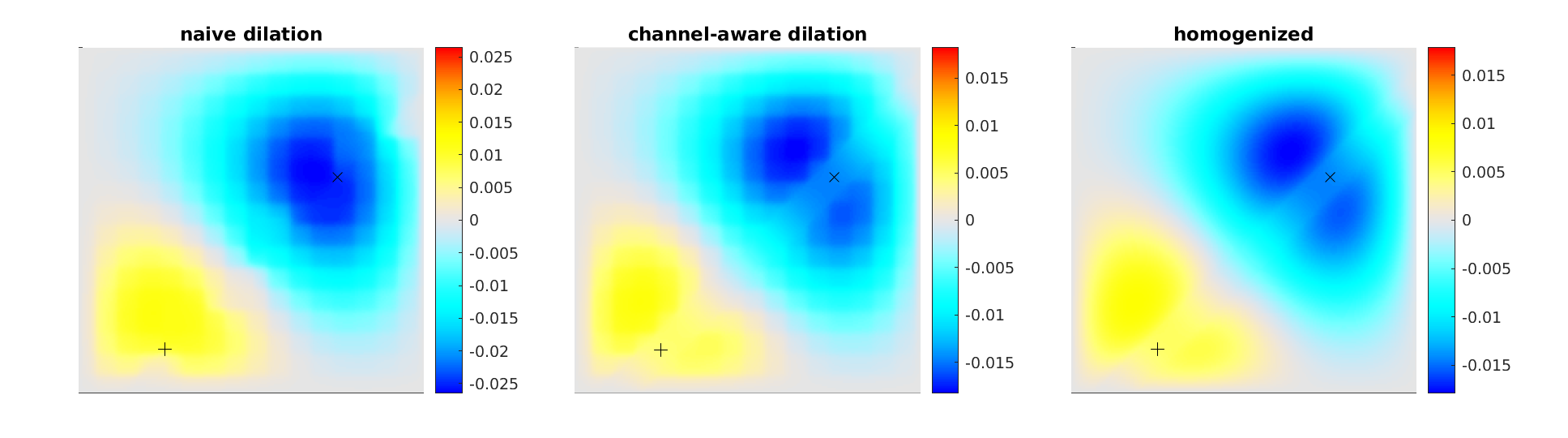}
\par\end{centering}
\label{fig:channel aware-dilation solution}}
\par\end{centering}
\begin{centering}
\subfloat[Comparing the velocity field $v_{1}^{\epsilon}:=e_{1}\cdot A^{\epsilon}\nabla u^{\epsilon}$ from the naive/structure-aware dilated system with the homogenized solution $\overline{v}_{1}:=e_{1}\cdot\overline{A}\nabla u_{0}$. ]{\begin{centering}
\includegraphics[width=0.9\linewidth]{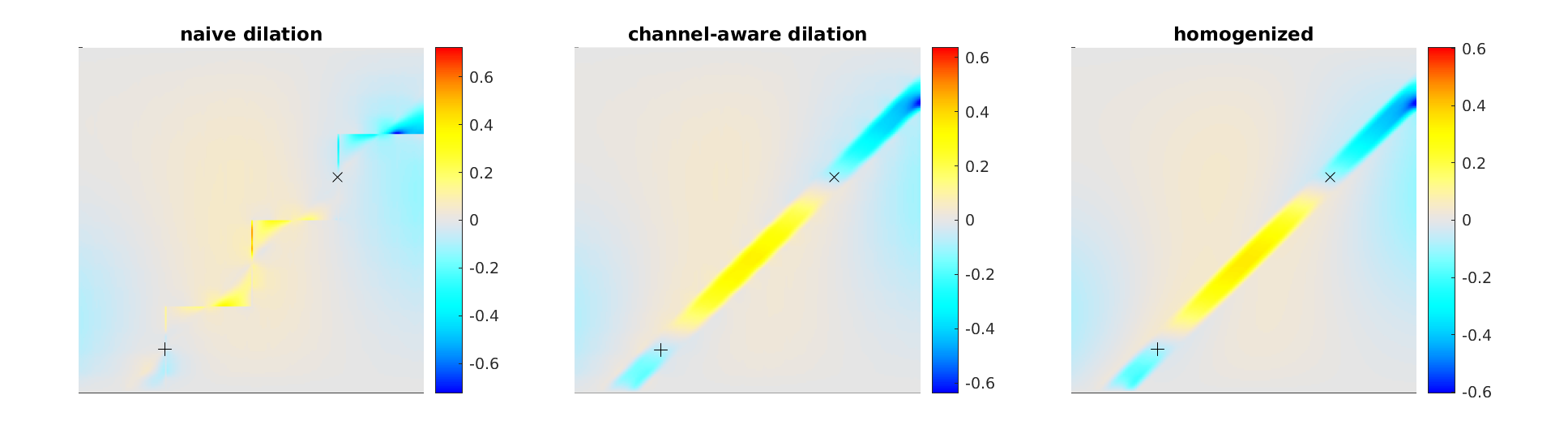}
\par\end{centering}
\label{fig:channel aware-dilation flux}}
\par\end{centering}
\caption{Illustration of the channel effect on coefficients and solutions.}

\label{fig:channel example}
\end{figure}

To conclude this section, we aim to address pitfalls that arise when the domain of interest may not exhibit the desired well-behaved heterogeneity. One example is the simulation of flow in porous media using Darcy's law, where the permeability might be drastically different in channels with known locations to the surrounding heterogeneous rock. Another example is heat conduction in composites containing isolated high-conductivity sheets or rods. These structures operate at a mesoscopic scale that is too large for the dilation operator to manipulate, thereby motivating the structure-aware dilation operator:
\begin{equation}
\cD_{L,m}^{\text{s}}\left(A^{\epsilon}\right):=A_{s}+\cD_{L,m}\left(A_{o}\right),\ \text{given}\,A^{\epsilon}\left(x\right)=A_{s}\left(x\right)+A_{o}\left(x,\frac{x}{\epsilon}\right)\label{eq:channel dilation operator}
\end{equation}
where
\begin{itemize}
\item the structure component $A_{s}\left(x\right)$ captures the location of abrupt changes in permeability or conductivity, and
\item the oscillation component $A_{o}\left(x,x/\epsilon\right)$ is in the form of 2-scale functions.
\end{itemize}
The decomposition (Eqn. \ref{eq:channel dilation operator}) can be acquired either through a priori knowledge or solely through a data-driven approach utilizing empirical mode decomposition. We display an example with a channel in Fig. \ref{fig:channel example}. The structure-aware dilation method prevents misalignments between the channel boundary and the mesoscale grid, thereby reducing numerical errors and enhancing the preservation of the flux passing through the structure.

\paragraph*{Numerical experiments}

\sepia{According to \cite[Theorem. 2.1]{allaireHomogenizationStokesFlow1989}, for Stokes equation} in the steady state, the velocity and pressure field 2-scale converges to the homogenized solution as the scale parameter $\epsilon$ vanishes. We are thus motivated to verify if the flux variable $v^{\epsilon}:=A^{\epsilon}\nabla u^{\epsilon}$ converges to the homogenized counterpart as $\epsilon\to0$. In particular, we aim to examine the scenario where there is enhanced permeability within the channel, motivated by physical principles (fluid movement facilitated within the channel, resulting in increased flux). The theoretical aspect of this convergence is left for further exploration.

The pressure field $u^{\epsilon}$ looks similar in both cases (Fig. \ref{fig:channel aware-dilation solution}), but the naive dilation fails to capture the nuisances inside the channel upon a closer examination. The discrepancy is fully revealed when we compute the velocity field (Fig. \ref{fig:channel aware-dilation flux}); the naive dilation method suffers from artificial continuities along the mesoscopic boundaries while the structure-aware method is able to recover the homogenized solution within small error.

\begin{figure}[tbph]
\begin{centering}
\subfloat[Convergence of relative error in solution $u$ w.r.t. $m$ (left) and $L$ (right). Notice that there is no big differences between the naive dilation and the structure-aware version.]{\begin{centering}
\includegraphics[width=0.9\linewidth]{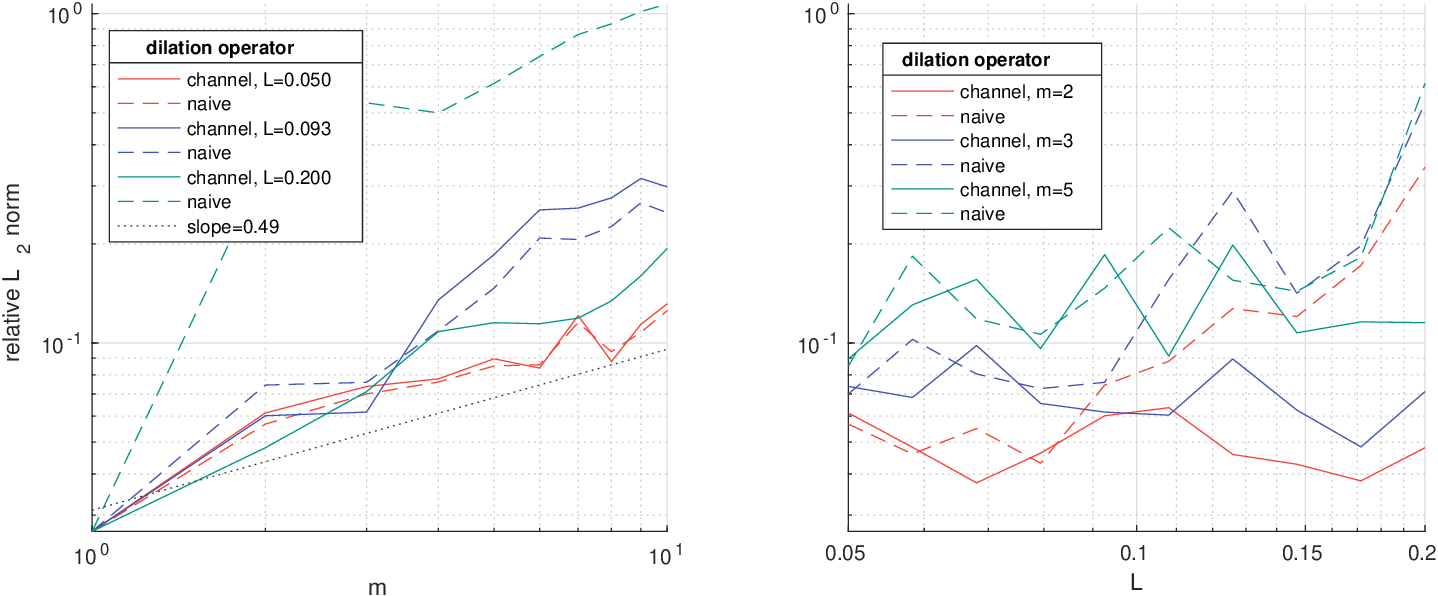}
\par\end{centering}

\label{fig:channel u convergence}}
\par\end{centering}
\begin{centering}
\subfloat[Convergence of relative error in flux $v_{1}^{\epsilon}=e_{1}\cdot A^{\epsilon}\nabla u^{\epsilon}$ w.r.t. $m$ (left) and $L$ (right). Notice that the naive dilation is worse than the structure-aware version.]{\begin{centering}
\includegraphics[width=0.9\linewidth]{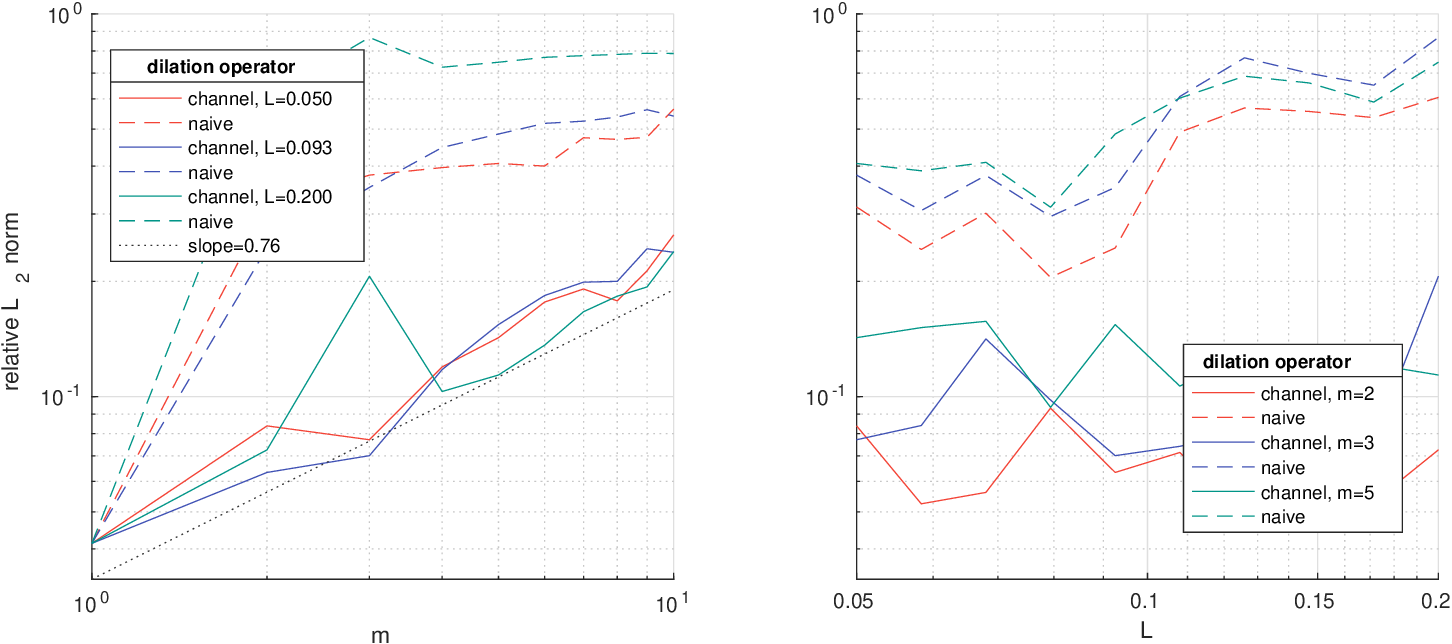}
\par\end{centering}

\label{fig:channel flux1 convergence}}
\par\end{centering}
\caption{Comparison of convergence behavior between naive and structure-aware dilation. The following system is tested: $A_{s}=\left(1+\eta_{c}\psi\left(\left(x_{2}-kx_{1}-b\right)/\sqrt{k^{2}+1};\epsilon_{c},s\right)\right)I_{2}$, $A_{o}=\eta_{o}\,\text{diag}\left(\sin\left(2\pi x_{1}/\epsilon\right),\sin\left(2\pi x_{2}/\epsilon\right)\right)$, $f=\sum\pm\exp\left(-\frac{1}{2}\times0.2^{2}\times\protect\abs{x-c_{\pm}}^{2}\right)$ with $\epsilon_{c}=0.03$, $k=1$, $b=-1/8$, $\eta_{c}=9$, $\eta_{o}=0.6$, $\epsilon=1/32$, $s=0.5$, $c_{+}=\left(\frac{1}{4},\frac{1}{4}k+b\right)$, $c_{-}=\left(\frac{3}{4},\frac{3}{4}k+b\right)$. The \textquotedblleft bridge function\textquotedblright{} $\psi\left(y;\epsilon,s\right)$ is defined as 0 for $\protect\abs{y/\epsilon}\le s$, $\left[\left(\protect\abs{y/\epsilon}-s\right)^{2}-1\right]^{2}$ for $s<\protect\abs{y/\epsilon}\le s+1$, and 1 otherwise. The coefficient field is visualized in Fig. \ref{fig:channel aware-dilation coeff}.}
\end{figure}

We continue to examine the overall error by the two dilation methods under different parameters of $m$ and $L$. As is shown in Fig. \ref{fig:channel u convergence}, these methods do not make a big difference in terms of the pressure field $u^{\epsilon}$. However, when it comes to the flux field $v^{\epsilon}$ (Fig. \ref{fig:channel flux1 convergence}), the structure-aware dilation achieves a lower error compared to naive dilation. The structure-aware dilation achieves a linear scaling w.r.t. to the scaling factor $m$ as in Sec. \ref{subsec:integrated-test}. The naive dilation, on the other hand, has a much larger error even when we take a scaling factor of $m=2$. The structure-aware dilation is also robust and insensitive to the $L$ parameter. An intuitive explanation for this difference is that the design of structure-aware dilation helps keep the channel intact so there are less artifacts around the channel.

\section{Discussion\label{sec:Discussion}}

We have proved a strong connection between the processes of averaging in dynamical systems and homogenization of elliptic operators in one dimension. This connection provides a foundation for extending the so-called Seamless method for oscillatory dynamical systems to multi-dimensional homogenization problems. The Seamless method manipulates the faster scale by introducing a mesoscale. It reduces the requirement of fine scale resolution, a requirement present in both oscillatory dynamical systems and homogenization, as well as in mathematical analysis and numerical methods. Here, we are considering efficient algorithms with sub-linear computational complexity relative to the inverse of the wavelength.

\sepia{We aim to discuss the relationships and differences between the three approaches mentioned in Sec. \ref{sec:Dilation-Solver}, focusing on whether scale identification is necessary and, if so, how the identification is performed.}
\begin{enumerate}
\item The first category, closely related to the Seamless method\cite{E2009}, requires an explicit separation between fast and slow equations so that it is possible to alter the effective $\epsilon$ parameter by manipulating the clock for the slow variables. Based on this, we propose the method of partial dilation in Sec. \ref{subsec:Averaging-and-partial-dilation} that increases the parameter $\epsilon$ in $A\left(x,x/\epsilon\right)$ to achieve a similar effect.
\item Methods in the second category share a lot of common ideas with these in the first category, except that the scale-separated form $A\left(x,x/\epsilon\right)$ is no longer available for the numerical procedure. Instead, the separation and identification of the slow and fast variables are conducted by tools from theories regarding signal processing. Then, the inferred form $\widetilde{A}\left(x,x/\epsilon\right)$ provides a proxy of the observations $A^{\epsilon}\left(x_{i}\right)$ where the scale parameter $\epsilon$ can be adjusted for further purposes. The details are elaborated in Sec. \ref{subsec:Numerical-inversion} and \ref{subsec:integrated-test}.
\item The third category does not assume a scale-separated form either; however, unlike the second category, these methods do not aim to recover $A\left(x,x/\epsilon\right)$. In contrast to the exact dilation effect as in the partial dilation method, an approximation is built for the coefficient via local spatial dilation. The approximated operator corresponds to the operator yielded by partial dilation up to the local dilation error.
\end{enumerate}
We shall point out that the local dilation method resembles FLAVORS in how the two different scales interact with each other. The FLAVORS integrator alternates between two stages that uses different time scales, one significantly smaller than $\epsilon$ (called ``transient phase'' for simplicity) and the other one larger than $\epsilon$ (``drifting phase''). The fast variable is only updated in the transient phase but held still in the drifting phase. The alternating pattern between short and large steps leads to a numerical homogenization effect. In local dilation, the transform on coefficients is equivalent to the transform on the coordinates, leading to an effective mesh structure that alternates between short and long intervals. This provides an another angle of understanding the local dilation solver. While this approach is always applicable, it is less accurate than the other two techniques and the error analysis is more complex. The algorithms are supported by illuminating numerical examples.

There are natural continuations of the research presented in this paper. In particular, the second and third methodologies mentioned above may be explored in connection to realistic applications. The third technique could be enhanced through regularization and adaptive meshing to mitigate the impact of discontinuities in its current version. It would also be interesting with a more detail analysis. \sepia{Furthermore, the partial/local dilation operator can be applied to elliptic operators of higher orders with minor modifications, provided the following conditions are met: a quantitative homogenization error (Eqn. \ref{eq:dilated-system-homogenization-error}) must be established, and the commutative condition (Lem. \ref{lem:interchangability}) must be verified.}

\bibliographystyle{siam}
\bibliography{zotero_ref}

\end{document}

%% file: math_shorthand.tex
\begin{comment}
mathbb, mathcal
\end{comment}

\global\long\def\bB{\mathbb{B}}%

\global\long\def\bC{\mathbb{C}}%

\global\long\def\bE{\mathbb{E}}%

\global\long\def\bF{\mathbb{F}}%

\global\long\def\bK{\mathbb{K}}%

\global\long\def\bN{\mathbb{N}}%

\global\long\def\bP{\mathbb{P}}%

\global\long\def\bQ{\mathbb{Q}}%

\global\long\def\bR{\mathbb{R}}%

\global\long\def\bS{\mathbb{S}}%

\global\long\def\bT{\mathbb{T}}%

\global\long\def\bZ{\mathbb{Z}}%

\global\long\def\cA{\mathcal{A}}%

\global\long\def\cB{\mathcal{B}}%

\global\long\def\cC{\mathcal{C}}%

\global\long\def\cD{\mathcal{D}}%

\global\long\def\cE{\mathcal{E}}%

\global\long\def\cF{\mathcal{F}}%

\global\long\def\cG{\mathcal{G}}%

\global\long\def\cH{\mathcal{H}}%

\global\long\def\cI{\mathcal{I}}%

\global\long\def\cJ{\mathcal{J}}%

\global\long\def\cK{\mathcal{K}}%

\global\long\def\cL{\mathcal{L}}%

\global\long\def\cLp#1{\mathcal{L}^{#1}}%

\global\long\def\cLpp#1{\mathcal{L}_{+}^{#1}}%

\global\long\def\cLsimp{\mathcal{L}_{simp}^{0}}%

\global\long\def\cM{\mathcal{M}}%

\global\long\def\cN{\mathcal{N}}%

\global\long\def\cO{\mathcal{O}}%

\global\long\def\cP{\mathcal{P}}%

\global\long\def\cR{\mathcal{R}}%

\global\long\def\cS{\mathcal{S}}%

\global\long\def\cT{\mathcal{T}}%

\global\long\def\cU{\mathcal{U}}%

\global\long\def\cV{\mathcal{V}}%

\global\long\def\cW{\mathcal{W}}%

\global\long\def\cY{\mathcal{Y}}%

\global\long\def\cZ{\mathcal{Z}}%

\global\long\def\fd{\mathscr{\mathfrak{d}}}%

\global\long\def\fm{\mathscr{\mathfrak{m}}}%

\global\long\def\fp{\mathscr{\mathfrak{p}}}%

\global\long\def\fq{\mathscr{\mathfrak{q}}}%

\global\long\def\fu{\mathfrak{u}}%

\global\long\def\fv{\mathfrak{v}}%

\global\long\def\fH{\mathscr{\mathfrak{H}}}%

\global\long\def\fL{\mathscr{\mathfrak{L}}}%

\global\long\def\fK{\mathscr{\mathfrak{K}}}%

\global\long\def\fN{\mathscr{\mathfrak{N}}}%

\global\long\def\fS{\mathscr{\mathfrak{S}}}%

\global\long\def\sD{\mathscr{D}}%

\global\long\def\sH{\mathscr{H}}%

\global\long\def\sK{\mathscr{K}}%

\global\long\def\sF{\mathscr{F}}%

\begin{comment}
norm, abs
\end{comment}

\global\long\def\norm#1{\left\Vert #1\right\Vert }%

\global\long\def\np#1#2{\left\Vert #1\right\Vert _{#2}}%

\global\long\def\nlp#1#2{\left\Vert #1\right\Vert _{L^{#2}}}%

\global\long\def\abs#1{\left|#1\right|}%

\global\long\def\inv#1{#1^{-1}}%

\begin{comment}
Linear algebra
\end{comment}

\global\long\def\adjoint#1{#1^{*}}%

\global\long\def\annihilator#1{#1^{\circ}}%

\global\long\def\annihilatee#1{#1^{\perp}}%

\global\long\def\unaryop#1{#1\left(\cdot\right)}%

\global\long\def\binaryop#1{#1\left(\cdot,\cdot\right)}%

\begin{comment}
other
\end{comment}

\global\long\def\comp#1#2{#1\circ#2}%

\global\long\def\converge#1{\overset{#1}{\joinrel\longrightarrow}}%

\global\long\def\define{\triangleq}%

\global\long\def\enum#1#2{\left\{  #1_{1},\dots,#1_{#2}\right\}  }%

\global\long\def\enumvec#1#2{\left(#1_{1},\dots,#1_{#2}\right)}%

\global\long\def\enuminf#1{\left\{  #1_{1},#1_{2}\dots\right\}  }%

\global\long\def\equivalent{\Longleftrightarrow}%

\global\long\def\substitute#1{\overset{#1}{\joinrel===}}%

\global\long\def\tensor{\otimes}%

\begin{comment}
rm-font operators
\end{comment}

\global\long\def\liminf#1{\underset{#1}{\operatorname{lim\,inf}}}%

\global\long\def\limsup#1{\underset{#1}{\operatorname{lim\,sup}}}%

\global\long\def\essinf#1{\underset{#1}{\operatorname{ess\,inf}}}%

\global\long\def\esssup#1{\underset{#1}{\operatorname{ess\,sup}}}%

\global\long\def\sgn{\operatorname{sgn}}%

\global\long\def\spanset{\operatorname{span}}%

\global\long\def\transpose#1{{#1}^{T}}%

\global\long\def\Null{\operatorname{Null}}%

\global\long\def\Range{\operatorname{Range}}%

\global\long\def\pv{\operatorname{p.v.}}%

\global\long\def\tr{\operatorname{\mathbf{Tr}}}%

\global\long\def\io{\operatorname{i.o.}}%

\global\long\def\ae{\operatorname{a.e.}}%

\global\long\def\as{\operatorname{a.s.}}%

\global\long\def\d#1{\operatorname{d}#1}%

\global\long\def\D#1{\operatorname{D}#1}%

\global\long\def\Db#1{\operatorname{D}\left[#1\right]}%

\global\long\def\cov{\operatorname{cov}}%

\global\long\def\supp{\operatorname{supp}}%